\documentclass[9pt]{article}
\usepackage{mathrsfs}
\usepackage{amsthm}
\usepackage{amssymb}
\usepackage{amsmath}
\usepackage{graphicx}
\usepackage{color}
\usepackage{amsfonts}
\usepackage{float}
\usepackage{cite}
\usepackage [latin1]{inputenc}
\usepackage[text={152mm,228mm},left=38mm,vmarginratio=1:1]{geometry}
\newtheorem{theorem}{Theorem}[section]
\newtheorem{remark}[theorem]{Remark}
\newtheorem{lemma}[theorem]{Lemma}

\numberwithin{equation}{section}
\normalsize

\begin{document}
\title{\textbf{Central limit theorems for additive functionals of long-range zero-range processes}}

\author{Xiaofeng Xue \thanks{\textbf{E-mail}: xfxue@bjtu.edu.cn \textbf{Address}: School of Mathematics and Statistics, Beijing Jiaotong University, Beijing 100044, China.}\\ Beijing Jiaotong University}

\date{}
\maketitle

\noindent {\bf Abstract:} In this paper, we extend the central limit theorem of the additive functional of the nearest-neighbor zero-range process given in \cite{Quastel2002} to the long-range case. Our main results show that in several cases the limit processes are driven by fractional Brownian motions with Hurst parameters in $(1/2, 3/4]$. A local central limit theorem of the long-range random walk and a relaxation to equilibrium theorem of the long-range zero-range process play the key roles in the proofs of our main results.

\quad


\section{Introduction}\label{section one}
In this paper, we are concerned with the long-range zero-range processes on lattices. We aim to extend the central limit theorems for the additive functionals of the nearest-neighbor zero-range processes given in \cite{Quastel2002} to the long-range case. We first give some notations for later use. For integer $d\geq 1$ and any vertex $x$ on the $d$-dimension lattice $\mathbb{Z}^d$, we denote by $\|x\|_2$ the $l_2$-norm of $x$, i.e.,
\[
\|x\|_2=\sqrt{\sum_{j=1}^dx_j^2}
\]
for any $x=(x_1, \ldots, x_d)\in \mathbb{Z}^d$. We denote by $\mathbb{N}$ the set of non-negative integers, i.e., $\mathbb{N}=\{0, 1, 2, \ldots\}$. Throughout this paper, we assume that $c: \mathbb{N}\rightarrow [0, +\infty)$ is a given function such that $c(0)=0$ and
\begin{equation}\label{equ 1.1 basic assumption of c}
0<\inf_{k\in \mathbb{N}}\left(c(k+1)-c(k)\right)\leq \sup_{k\in \mathbb{N}}\left(c(k+1)-c(k)\right)<+\infty.
\end{equation}
We denote by $\mathbf{0}$ the origin of $\mathbb{Z}^d$, i.e., $\mathbf{0}=(0, 0, \ldots, 0)$.

Now we recall the definition of the long-range zero-range process. For given integer $d\geq 1$ and constant $\alpha>0$, the long-range zero-range process $\{\eta_t\}_{t\geq 0}$ on $\mathbb{Z}^d$ with parameter $\alpha$ is a continuous-time Markov process with state space $\mathbb{N}^{\mathbb{Z}^d}$. The generator $\mathcal{L}$ of $\{\eta_t\}_{t\geq 0}$ is given by
\begin{equation}\label{equ 1.2 generator}
\mathcal{L}f(\eta)=\sum_{x\in \mathbb{Z}^d}\sum_{y\neq x}c(\eta(x))\|y-x\|_2^{-(d+\alpha)}\left(f(\eta^{x, y})-f(\eta)\right)
\end{equation}
for any $\eta\in \mathbb{N}^{\mathbb{Z}^d}$ and local $f: \mathbb{N}^{\mathbb{Z}^d}\rightarrow \mathbb{R}$, where
\[
\eta^{x, y}(z)=
\begin{cases}
\eta(z) & \text{~if~}z\neq x, y,\\
\eta(x)-1 & \text{~if~}z=x,\\
\eta(y)+1 & \text{~if~}z=y.
\end{cases}
\]
When we need to emphasize the dependence of $\mathcal{L}$ on $d, \alpha$, we write $\mathcal{L}$ as $\mathcal{L}^{d, \alpha}$. According to the definition of $\mathcal{L}$, the long-range zero-range process evolves as follows. There are several particles on each vertex on $\mathbb{Z}^d$. For any $x\in \mathbb{Z}^d$, each particle on $x$ jumps to $y\in \mathbb{Z}^d$ at rate $\frac{c(\eta(x))}{\eta(x)}\|y-x\|_2^{-(d+\alpha)}$, where $\eta(x)$ is the number of particles on $x$. Note that, by Assumption \eqref{equ 1.1 basic assumption of c}, there exists $c^+<+\infty$ such that $c(k)\leq kc^+$ for all $k\in \mathbb{N}$ and then
each particle jumps at rate at most
\[
c^+\sum_{y\neq 0}\|y\|_2^{-(d+\alpha)}<+\infty.
\]
Therefore, the process $\{\eta_t\}_{t\geq 0}$ is well-defined. When $c(k)=kc(1)$ for all $k\in \mathbb{N}$, the model reduces to the long-range Ehrenfest model on $\mathbb{Z}^d$, where all particles perform independent long-range random walks on $\mathbb{Z}^d$.

For later use, for any probability distribution $\mu$ on $\mathbb{N}^{\mathbb{Z}^d}$, we denote by $\mathbb{P}_\mu$ the probability measure of $\{\eta_t\}_{t\geq 0}$ starting from $\mu$. We further denote by $\mathbb{E}_\mu$ the expectation with respect to $\mathbb{P}_\mu$. When $\mu$ is the Dirac measure concentrated on $\eta\in \mathbb{N}^{\mathbb{Z}^d}$, we write $\mathbb{P}_\mu$ and $\mathbb{E}_\mu$ as $\mathbb{P}_\eta$ and $\mathbb{E}_\eta$ respectively.

Now we recall a family of reversible distributions of $\{\eta_t\}_{t\geq 0}$. For given parameter $\beta>0$, we denote by $\nu$ the probability measure on $\mathbb{N}^{\mathbb{Z}^d}$ such that $\{\eta(x)\}_{x\in \mathbb{Z}^d}$ are independent under $\nu$ and
\[
\nu\left(\eta:~\eta(x)=k\right)=\frac{1}{Z(\beta)}\frac{\beta^k}{c(k)!}
\]
for any $x\in \mathbb{Z}^d$ and $k\in \mathbb{N}$, where
\[
c(k)!=
\begin{cases}
\prod_{l=1}^kc(l) & \text{~if~}k\geq 1,\\
1 & \text{~if~} k=0
\end{cases}
\]
and $Z(\beta)=\sum_{k=0}^{+\infty}\frac{\beta^k}{c(k)!}$. Note that, by Assumption \eqref{equ 1.1 basic assumption of c}, there exists $c^-\in (0, +\infty)$ such that $c(k)\geq kc^-$ for all $k\in \mathbb{N}$ and hence $Z(\beta)\leq e^{\beta/{c^-}}<+\infty$.

According to the definition of $\nu$, for any different $x, y\in \mathbb{Z}^d$ and $k\geq 1, l\geq 0$, we have
\[
\nu\left(\eta:~\eta(x)=k, \eta(y)=l\right)c(k)\|x-y\|_2^{-(d+\alpha)}=\nu\left(\eta:~\eta(x)=k-1, \eta(y)=l+1\right)c(l+1)\|y-x\|_2^{-(d+\alpha)}
\]
and hence $\nu$ is a reversible distribution of $\{\eta_t\}_{t\geq 0}$. We denote by $\gamma=\gamma(\beta)$ the mean of $\eta(\mathbf{0})$ under $\nu$, i.e.,
\[
\gamma=\mathbb{E}_\nu\eta(\mathbf{0}).
\]
It is obvious that $\gamma$ is strictly increasing in $\beta$. Therefore, for any $\gamma>0$, there exists a unique $\beta=\beta(\gamma)$ to make $\gamma=\gamma(\beta)$. It is common to parameterize $\nu$ by $\gamma$ instead of $\beta$ and hence we write $\nu$ as $\nu_\gamma$ when $\nu$ satisfies that $\mathbb{E}_\nu\eta(\mathbf{0})=\gamma$. According to the definition of $\nu_\gamma$, it is not difficult to show that
\begin{equation}\label{equ 1.3}
\mathbb{E}_{\nu_\gamma}c(\eta(\mathbf{0}))=\beta=\beta(\gamma)
\end{equation}
and
\begin{equation}\label{equ 1.4}
{\rm Var}_{\nu_\gamma}\left(\eta(\mathbf{0})\right)=\frac{\beta(\gamma)}{\frac{d}{d\gamma}\beta(\gamma)}.
\end{equation}
For detailed checks of \eqref{equ 1.3} and \eqref{equ 1.4}, see Section 5 of \cite{kipnis+landim99}.

For a local function $V: \mathbb{N}^{\mathbb{Z}^d}\rightarrow \mathbb{R}$ such that $\mathbb{E}_{\nu_\gamma}V=0$, the process $\left\{\int_0^tV(\eta_s)ds\right\}_{t\geq 0}$ is called an additive functional of $\{\eta_t\}_{t\geq 0}$. A special case of the additive functional is the occupation time process on $\mathbf{0}$, where $V(\eta)=\eta(\mathbf{0})-\gamma$. Reference \cite{Quastel2002} gives central limit theorems of the additive functionals of the nearest-neighbor zero-range process, where each particle on vertex $x$ jumps to each neighbor $y$ of $x$ at rate $c(\eta(x))/\eta(x)$. By Theorem 1.2 of \cite{Quastel2002}, for the nearest-neighbor zero-range process $\{\xi_t\}_{t\geq 0}$ on $\mathbb{Z}^d$ starting from $\nu_\gamma$, the rescaled
additive functional $\left\{\frac{1}{A_d(N)}\int_0^{tN}V(\xi_s)ds\right\}_{0\leq t\leq T}$ converges weakly, with respect to the uniform topology, to a Gaussian process $\{X_t^d\}_{0\leq t\leq T}$ as $N\rightarrow+\infty$, where
\[
A_d(N)=
\begin{cases}
N^{\frac{3}{4}} & \text{~if~} d=1,\\
\sqrt{N\log N} & \text{~if~} d=2,\\
\sqrt{N} & \text{~if~} d\geq 3
\end{cases}
\]
and $\{X_t^d\}$ is a Brownian motion when $d\geq 2$ or a fractional Brownian motion with Hurst parameter $3/4$ when $d=1$. The main theorem of this paper extends the above central limit theorem to the long-range case. For mathematical details, see Section \ref{section two}.

We are also inspired by \cite{Bernardin2016}, where the central limit theorem of the occupation time of the long-range symmetric exclusion process (SEP) is given. For the long-range SEP on $\mathbb{Z}^d$ with parameter $\alpha>0$, the particle on vertex $x$ jumps to each vacant vertex $y$ at rate $\|y-x\|_2^{-(d+\alpha)}$ and any jump to a occupied vertex is suppressed to ensure that there is at most one particle on each vertex. In the main theorem (Theorem 2.11) of \cite{Bernardin2016}, the most interesting part is that, for the long-range SEP on $\mathbb{Z}^1$ with parameter $\alpha\in (1, 2)$, the weak limit of the rescaled occupation time is driven by the fractional Brownian motion with Hurst parameter $1-1/2\alpha$. In this paper, we show that the additive functional of the long-range zero-range process on $\mathbb{Z}^1$ with $\alpha\in (1, 2)$ has similar asymptotic behavior. For mathematical details, see Section \ref{section two}.

Central limit theorems of additive functionals of long-range interacting particle systems are also discussed for the voter model. For the long-range voter model on $\mathbb{Z}^d$ with parameter $\alpha>0$, each vertex is a supporter of two opposite opinions and a vertex $x$ adopts another vertex $y$'s opinion at rate $\|x-y\|_2^{-(d+\alpha)}$. Reference \cite{Xue2025} shows that the central limit theorem of the occupation time of the long-range voter model performs phase transitions divided into eight different cases. For mathematical details, see Theorems 2.1-2.4 of \cite{Xue2025}.

\section{Main results} \label{section two}
In this paper, we give our main results. As a preliminary, we first introduce some definitions and notations. For integer $d\geq 1$ and $\alpha>0$, we denote by $\{\hat{X}_t^{d, \alpha}\}_{t\geq 0}$ the stable process with state space $\mathbb{R}^d$ and generator $\hat{\mathcal{G}}^{d, \alpha}$ defined as follows. When $\alpha\in (0, 2)$,
\[
\hat{\mathcal{G}}^{d, \alpha}f(u)=\int_{\mathbb{R}^d}\frac{f(u+v)-f(u)}{\|v\|_2^{d+\alpha}}dv
\]
for any $u\in \mathbb{R}^d$ and $f\in C_c^2(\mathbb{R}^d)$. When $\alpha>2$,
\[
\hat{\mathcal{G}}^{d, \alpha}f(u)=\frac{1}{2}\sum_{i=1}^d\sum_{j=1}^d\left(\sum_{y\in \mathbb{Z}^d\setminus\{\mathbf{0}\}}\frac{y_iy_j}{\|y\|_2^{d+\alpha}}\right)\partial_{u_iu_j}^2f(u)
\]
for any $u\in \mathbb{R}^d$ and $f\in C_c^2(\mathbb{R}^d)$, where $u_i$ is the $i$th coordinate of $u$. When $\alpha=2$,
\[
\hat{\mathcal{G}}^{d, 2}f(u)=\frac{1}{2}\sum_{i=1}^d\sum_{j=1}^dK_{d,2}^{i,j}\partial_{u_iu_j}^2f(u)
\]
for any $u\in \mathbb{R}^d$ and $f\in C_c^2(\mathbb{R}^d)$, where
\[
K_{d,2}^{i,j}=\frac{1}{2}\int_{\{v\in \mathbb{R}^d:~\|v\|_2=1\}}v_iv_jdS
\]
and $dS$ is the surface integral. We denote by $f_t^{d, \alpha}$ the probability density of $\hat{X}_t^{d, \alpha}$ conditioned on $\hat{X}_0^{d, \alpha}=\mathbf{0}$. Note that, according to the self-similarity of the stable process (see Section 3.7 of \cite{Dur2010}), $f_t^{d, \alpha}(\mathbf{0})=t^{-\frac{d}{2}}f_1^{d, \alpha}(\mathbf{0})$ when $\alpha\geq 2$ and
\[
f_t^{d, \alpha}(\mathbf{0})=t^{-\frac{d}{\alpha}}f_1^{d, \alpha}(\mathbf{0}).
\]
when $\alpha\in (0, 2)$.

For $\theta\in (0, 1)$, we denote by $\{B_t^\theta\}_{t\geq 0}$ the fractional Brownian motion with Hurst parameter $\theta$, i.e., $\{B_t^\theta\}_{t\geq 0}$ is a Gaussian process with continuous sample path, mean zero and covariance function given by
\[
{\rm Cov}\left(B_t^\theta, B_s^\theta\right)=\frac{1}{2}\left(t^{2\theta}+s^{2\theta}-|t-s|^{2\theta}\right)
\]
for any $t,s>0$. We further denote $B_t^{\frac{1}{2}}$ by $\mathcal{W}_t$, i.e., $\{\mathcal{W}_t\}_{t\geq 0}$ is the standard Brownian motion starting from $0$.

For any $V: \mathbb{N}^{\mathbb{Z}^d}\rightarrow \mathbb{R}$, if there exists $M=M_V<+\infty$ such that $V$ only depends on
\[
\left\{\eta(x):~\|x\|_2\leq M\right\},
\]
then we call $V$ a local function. Furthermore, for any local $V$, if there exists $k\in (0, +\infty)$ such that
\[
|V(\eta)|\leq \left(\sum_{x:\|x\|_2\leq M_V}\eta(x)\right)^k
\]
for any $\eta\in \mathbb{N}^{\mathbb{Z}^d}$, then we call that $V$ is with polynomial bound.

For any $\gamma>0$ and local $V: \mathbb{N}^{\mathbb{Z}^d}\rightarrow \mathbb{R}$ with polynomial bound, we define
\[
\overline{V}(\gamma)=\mathbb{E}_{\nu_\gamma}V.
\]
For simplicity, we denote $\frac{d}{d\tilde{\gamma}}\overline{V}(\tilde{\gamma})\Big|_{\tilde{\gamma}=\gamma}$ by $\overline{V}^\prime(\gamma)$. For example, when $V(\eta)=\eta(\mathbf{0})$, we have $\overline{V}(\gamma)=\gamma$ and $\overline{V}^\prime(\gamma)=1$.

For $\gamma>0$ and local $V:\mathbb{N}^{\mathbb{Z}^d}\rightarrow\mathbb{R}$ such that $\overline{V}(\gamma)=0$, we denote by $\sigma_\gamma(V)$ the term
\[
\sqrt{2\int_0^{+\infty}\mathbb{E}_{\nu_\gamma}\left(V(\eta_0)V(\eta_s)\right)ds}.
\]
By a crucial variance estimation given later (see Lemma \ref{lemma relaxation to equibrium}) and the reversibility of $\nu_\gamma$, we have that
\begin{equation}\label{equ 2.1 order of variance}
{\rm Var}_{\nu_\gamma}\left(\mathbb{E}_{\eta_0}V(\eta_t)\right)=\mathbb{E}_{\nu_\gamma}\left(V(\eta_0)V(\eta_{2t})\right)=O\left((h_\alpha(t))^{-d}\right)
\end{equation}
as $t\rightarrow+\infty$ for any local $V$ with polynomial bound satisfying $\overline{V}(\gamma)=0$, where
\[
h_\alpha(t)=
\begin{cases}
t^{1/\alpha} & \text{~if~} \alpha\in (0, 2),\\
\sqrt{t\log t} & \text{~if~} \alpha=2,\\
\sqrt{t} & {\text{~if~}} \alpha>2.
\end{cases}
\]

Therefore, when
\[
(d, \alpha)\in \{1\}\times (0, 1)\bigcup \{2\}\times (0, 2)\bigcup \{3, 4, \ldots\}\times (0, +\infty),
\]
$\sigma_\gamma(V)<+\infty$ for any local $V$ with polynomial bound satisfying $\overline{V}(\gamma)=0$.

Now we give our main result in the case $d=1$.

\begin{theorem}\label{theorem 2.1 main theorem d=1}
Let $T>0$, $\gamma>0$, local $V: \mathbb{N}^{\mathbb{Z}^1}\rightarrow \mathbb{R}$ with polynomial bound satisfy that $\overline{V}(\gamma)=0, \overline{V}^\prime(\gamma)\neq 0$ and $\{\eta_t\}_{t\geq 0}$ on $\mathbb{Z}^1$ with parameter $\alpha>0$ start from $\nu_\gamma$. If $\alpha\neq 1$, then
\[
\left\{\frac{1}{\Lambda_{1, \alpha}(N)}\int_0^{tN}V(\eta_s)ds:~0\leq t\leq T\right\}
\]
converges weakly, with respect to the uniform topology, to $\left\{Y_t^{1, \alpha}\right\}_{0\leq t\leq T}$ as $N\rightarrow+\infty$, where
\[
Y_t^{1, \alpha}=
\begin{cases}
\sigma_\gamma(V)\mathcal{W}_t & \text{~if\quad} \alpha<1,\\
\sqrt{\frac{2\alpha^2}{(\alpha-1)(2\alpha-1)}\frac{\left(\overline{V}^\prime(\gamma)\right)^2\beta(\gamma)}{\left(\beta^\prime(\gamma)\right)^{1+\frac{1}{\alpha}}}f_1^{1, \alpha}(\mathbf{0})}B_t^{1-\frac{1}{2\alpha}} & \text{~if\quad}1<\alpha<2,\\
\sqrt{\frac{8}{3}\frac{\left(\overline{V}^\prime(\gamma)\right)^2\beta(\gamma)}{\left(\beta^\prime(\gamma)\right)^{\frac{3}{2}}}f_1^{1, \alpha}(\mathbf{0})}B_t^{\frac{3}{4}} & \text{~if\quad}\alpha\geq 2
\end{cases}
\]
and
\[
\Lambda_{1, \alpha}(N)=
\begin{cases}
\sqrt{N} & \text{~if\quad} \alpha<1,\\
N^{1-\frac{1}{2\alpha}} & \text{~if\quad}1<\alpha<2,\\
N^{\frac{3}{4}}(\log N)^{-\frac{1}{4}} & \text{~if\quad}\alpha=2,\\
N^{\frac{3}{4}} & \text{~if\quad}\alpha>2.
\end{cases}
\]
If $\alpha=1$, then for any $0<t_1<t_2<\ldots<t_m$,
\[
\frac{1}{\sqrt{N\log N}}\left(\int_0^{t_1N}V(\eta_s)ds, \int_0^{t_2N}V(\eta_s)ds, \ldots, \int_0^{t_mN}V(\eta_s)ds\right)
\]
converges weakly to $\left(Y_{t_1}^{1, 1}, Y_{t_2}^{1, 1}, \ldots, Y_{t_m}^{1, 1}\right)$ as $N\rightarrow+\infty$, where
\[
Y_t^{1, 1}=\sqrt{2f_1^{1,1}(\mathbf{0})\left(\frac{\overline{V}^\prime(\gamma)}{\beta^\prime(\gamma)}\right)^2\beta(\gamma)}\mathcal{W}_t.
\]
\end{theorem}

We believe that, when $\alpha=1$, $\left\{\frac{1}{\sqrt{N\log N}}\int_0^{tN}V(\eta_s)ds:~0\leq t\leq T\right\}$ also converges weakly to
$\{Y_t^{1, 1}\}_{0\leq t\leq T}$ under the uniform topology. However, according to some technical reasons (see Remark \ref{remark 4.10}), currently we can only check the tightness of
\[
\left\{\frac{1}{\Lambda_{1, \alpha}(N)}\int_0^{tN}V(\eta_s)ds:~0\leq t\leq T\right\}_{N\geq 1}
\]
under the uniform topology in the case $\alpha\neq 1$. So, in the case $\alpha=1$, we only obtain a finite-dimensional marginal distribution CLT.

By \eqref{equ 2.1 order of variance}, when $d=1$,
\[
{\rm Var}_{\nu_\gamma}\left(\int_0^{tN}V(\eta_s)ds\right)=
\begin{cases}
O(1)tN & \text{~if~}\alpha<1,\\
O(1)tN\log N & \text{~if~} \alpha=1,\\
O(1)t^{2-\frac{1}{\alpha}}N^{2-\frac{1}{\alpha}} & \text{~if~} 1<\alpha<2,\\
O(1)t^{3/2}N^{3/2}(\log N)^{-1/2} & \text{~if~} \alpha=2,\\
O(1)t^{3/2}N^{3/2} & \text{~if~} \alpha>2
\end{cases}
\]
as $N\rightarrow+\infty$, which intuitively explains why the phase transitions in Theorem \ref{theorem 2.1 main theorem d=1} are divided into the aforesaid five parts.

Now we give our main theorem in the case $d=2$.

\begin{theorem}\label{theorem 2.2 main theorem d=2}
Let $T>0$, $\gamma>0$, local $V: \mathbb{N}^{\mathbb{Z}^2}\rightarrow \mathbb{R}$ with polynomial bound satisfy that $\overline{V}(\gamma)=0, \overline{V}^\prime(\gamma)\neq 0$ and $\{\eta_t\}_{t\geq 0}$ on $\mathbb{Z}^2$ with parameter $\alpha>0$ start from $\nu_\gamma$. If $\alpha<2$, then
\[
\left\{\frac{1}{\sqrt{N}}\int_0^{tN}V(\eta_s)ds:~0\leq t\leq T\right\}
\]
converges weakly, with respect to the uniform topology, to $\left\{\sigma_\gamma(V)\mathcal{W}_t\right\}_{0\leq t\leq T}$ as $N\rightarrow+\infty$.
If $\alpha\geq 2$, then for any $0<t_1<t_2<\ldots<t_m$,
\[
\frac{1}{\Lambda_{2, \alpha}(N)}\left(\int_0^{t_1N}V(\eta_s)ds, \int_0^{t_2N}V(\eta_s)ds, \ldots, \int_0^{t_mN}V(\eta_s)ds\right)
\]
converges weakly to $\left(Y_{t_1}^{2, \alpha}, Y_{t_2}^{2, \alpha}, \ldots, Y_{t_m}^{2, \alpha}\right)$ as $N\rightarrow+\infty$,
where
\[
Y_t^{2, \alpha}=
\sqrt{2f_1^{2,\alpha}(\mathbf{0})\left(\frac{\overline{V}^\prime(\gamma)}{\beta^\prime(\gamma)}\right)^2\beta(\gamma)}\mathcal{W}_t
\]
and
\[
\Lambda_{2, \alpha}(N)=
\begin{cases}
\sqrt{N\log (\log N)} & \text{~if\quad} \alpha=2,\\
\sqrt{N\log N} & \text{~if\quad} \alpha>2.
\end{cases}
\]
\end{theorem}

We believe that, when $\alpha\geq 2$, $\left\{\frac{1}{\Lambda_{2, \alpha}(N)}\int_0^{tN}V(\eta_s)ds:~0\leq t\leq T\right\}$ also converges weakly to $\{Y_t^{2, \alpha}\}_{0\leq t\leq T}$ under the uniform topology. Still, in this case we have not managed to check the tightness of $\left\{\frac{1}{\Lambda_{2, \alpha}(N)}\int_0^{tN}V(\eta_s)ds:~0\leq t\leq T\right\}_{N\geq 1}$ under the uniform topology and hence only obtain a finite-dimensional marginal distribution CLT.

By \eqref{equ 2.1 order of variance}, when $d=2$,
\[
{\rm Var}_{\nu_\gamma}\left(\int_0^{tN}V(\eta_s)ds\right)=
\begin{cases}
O(1)tN & \text{~if~}\alpha<2,\\
O(1)tN\log (\log N) & \text{~if~} \alpha=2,\\
O(1)tN\log N & \text{~if~} \alpha>2
\end{cases}
\]
as $N\rightarrow+\infty$, which intuitively explains why the phase transitions in Theorem \ref{theorem 2.2 main theorem d=2} are divided into the aforesaid three parts.

Now we give our main theorem in the case $d\geq 3$.
\begin{theorem}\label{theorem 2.3 main theorem d geq 3}
Let $T>0$,  $d\geq 3$, $\gamma>0$ and local $V: \mathbb{N}^{\mathbb{Z}^d}\rightarrow \mathbb{R}$ with polynomial bound satisfy that $\overline{V}(\gamma)=0, \overline{V}^\prime(\gamma)\neq 0$. If $\{\eta_t\}_{t\geq 0}$ on $\mathbb{Z}^d$ with parameter $\alpha>0$ starts from $\nu_\gamma$, then
\[
\left\{\frac{1}{\sqrt{N}}\int_0^{tN}V(\eta_s)ds:~0\leq t\leq T\right\}
\]
converges weakly, with respect to the uniform topology, to $\left\{\sigma_\gamma(V)\mathcal{W}_t\right\}_{0\leq t\leq T}$ as $N\rightarrow+\infty$.
\end{theorem}

By \eqref{equ 2.1 order of variance}, when $d\geq 3$,
\[
{\rm Var}_{\nu_\gamma}\left(\int_0^{tN}V(\eta_s)ds\right)=O(1)tN
\]
as $N\rightarrow+\infty$, which explains why there is no phase transition in Theorem \ref{theorem 2.3 main theorem d geq 3}.

The remainder of this paper is arranged as follows. In Section \ref{section three}, we give two crucial lemmas as preliminaries. The first lemma, which was proved in \cite{Xue2025}, is a local central limit theorem of the long-range random walk on $\mathbb{Z}^d$. The second one is an extension of the relaxation to equilibrium theorem of the nearest-neighbor zero-range process given in \cite{Janvresse1999} to the long-range case. The proofs of Theorems \ref{theorem 2.1 main theorem d=1} to \ref{theorem 2.3 main theorem d geq 3} are divided into Sections \ref{section four} and \ref{section five}, where we utilize the martingale decomposition strategy given in \cite{Kipnis1987} and the Poisson flow strategy given in \cite{Quastel2002}. By the martingale decomposition strategy, we can decompose the additive functional as a martingale plus a remainder term which only depends on the initial state. The CLT of the remainder term is easy to obtain since $\nu_\gamma$ is a product measure. By the Poisson flow strategy, we can write the martingale term as the Lebesgue-Stieltjes integral of a function with respect to a series of centered Poisson processes, where the function is related to the transition probabilities of the long-range random walk. Then, the CLT of the martingale term follows from the CLT of a random field generated by these centered Poisson processes and the aforesaid local central limit theorem of the long-range random walk. The CLT of the aforesaid random field follows from a routine quadratic variation computation, where the aforesaid relaxation to equilibrium theorem plays the key role. For mathematical details, see Sections \ref{section four} and \ref{section five}.

\section{Preliminaries} \label{section three}
In this section, we give two crucial lemmas for the proofs of Theorems \ref{theorem 2.1 main theorem d=1} to \ref{theorem 2.3 main theorem d geq 3}. We first introduce some notations for later use. For integer $d\geq 1$ and constant $\alpha>0$, we denote by $\{X_t^{d, \alpha}\}_{t\geq 0}$ the long-range random walk on $\mathbb{Z}^d$ with parameter $\alpha>0$, the generator $\mathcal{G}^{d, \alpha}$ of which is defined as
\[
\mathcal{G}^{d, \alpha}h(x)=\sum_{y\in \mathbb{Z}^d\setminus\{x\}}\|y-x\|_2^{-(d+\alpha)}\left(h(y)-h(x)\right)
\]
for any $x\in \mathbb{Z}^d$ and bounded $h: \mathbb{Z}^d\rightarrow \mathbb{R}$. We denote by $\{p_t^{d, \alpha}\}_{t\geq 0}$ the transition probabilities of $\{X_t^{d, \alpha}\}_{t\geq 0}$, i.e.,
\[
p_t^{d, \alpha}(x, y)=\mathbb{P}\left(X_t^{d, \alpha}=y\Big|X_0^{d, \alpha}=x\right)
\]
for any $x, y\in \mathbb{Z}^d$ and $t\geq 0$.

The following lemma, which was proved in \cite{Xue2025}, gives the local central limit theorem of $\{X_t^{d, \alpha}\}_{t\geq 0}$.
\begin{lemma}\label{lemma LCLT}
(\cite{Xue2025}). For any $t>0$,
\[
\lim_{s\rightarrow+\infty}\sup_{x\in \mathbb{Z}^d}\left|\left(h_\alpha(s)\right)^dp_{ts}^{d, \alpha}\left(\mathbf{0}, x\right)-f_t^{d, \alpha}\left(\frac{x}{h_\alpha(s)}\right)\right|=0,
\]
where $f_t^{d, \alpha}$ is the probability density of $\hat{X}_t^{d, \alpha}$ condition on $\hat{X}_0^{d, \alpha}=\textbf{0}$ defined as in Section \ref{section two} and
\[
h_{\alpha}(s)=
\begin{cases}
s^{\frac{1}{\alpha}} & \text{~if~}\alpha\in (0, 2),\\
\sqrt{s\log s} & \text{~if~} \alpha=2,\\
\sqrt{s} & \text{~if~} \alpha>2
\end{cases}
\]
for any $s\geq 0$.
\end{lemma}
Lemma \ref{lemma LCLT} extends the local central limit theorem of the simple random walk on $\mathbb{Z}^d$ (see Theorem 2.1.3 in Chapter 2 of \cite{Lawler2010} ) to the long-range case. For a rigorous proof of Lemma \ref{lemma LCLT}, see Section 3 of \cite{Xue2025}. Here we give an intuitive explanation of this lemma. It is not difficult to show that
\[
\lim_{s\rightarrow+\infty}\mathcal{G}_s^{d, \alpha}f(u)=\hat{\mathcal{G}}^{d, \alpha}f(u)
\]
for any $u\in \mathbb{R}^d$ and bounded $f: \mathbb{R}^d\rightarrow \mathbb{R}$, where
\[
\mathcal{G}_s^{d, \alpha}f(u)=s\sum_{y\in \mathbb{Z}^d\setminus\{\mathbf{0}\}}\|y\|_2^{-(d+\alpha)}\left(f\left(u+\frac{1}{h_\alpha(s)}y\right)-f(u)\right).
\]
Then, conditioned on $X_0^{d, \alpha}=\mathbf{0}$, $\frac{1}{h_\alpha(s)}X_{ts}^{d, \alpha}$ converges weakly, with respect to the Skorohod topology, to $\hat{X}_t^{d, \alpha}$ starting from $\mathbf{0}$ as $N\rightarrow+\infty$. Then, for sufficiently small $\epsilon>0$ and sufficiently large $s$,
\begin{align*}
(2\epsilon)^df_t^{d, \alpha}(u)&\approx \mathbb{P}\left(\hat{X}_t^{d, \alpha}\in \mathcal{B}(u, \epsilon)\Big|\hat{X}_0^{d, \alpha}=\mathbf{0}\right)
\\
&\approx\mathbb{P}\left(\frac{1}{h_\alpha(s)}X_{ts}^{d, \alpha}\in \mathcal{B}(u, \epsilon)\Big|X_0^{d, \alpha}=0\right)\\
&=\mathbb{P}\left(X_{ts}^{d, \alpha}\in \mathcal{B}(h_\alpha(s)u, h_\alpha(N)\epsilon)\Big|X_0^{d, \alpha}=0\right)\\
&\approx (2h_\alpha(s)\epsilon)^d\mathbb{P}\left(X_{ts}^{d, \alpha}=h_\alpha(s)u\Big|X_0^{d, \alpha}=0\right),
\end{align*}
where $\mathcal{B}(u, \epsilon)=\{v:\|v-u\|_\infty\leq \epsilon\}$. Hence, for large $s$,
\[
(h_\alpha(s))^dp_t^{d, \alpha}(\mathbf{0}, h_\alpha(s)u)\approx f_t^{d, \alpha}(u).
\]

Theorem 1.1 of \cite{Janvresse1999} gives a relaxation to equilibrium proposition of the nearest-neighbor zero-range process. Our next lemma extends this theorem to the long-range case.
\begin{lemma}\label{lemma relaxation to equibrium}
For any local $V: \mathbb{N}^{\mathbb{Z}^d}\rightarrow \mathbb{R}$ with polynomial bound satisfying that $\overline{V}(\gamma)=0$,
\[
\lim_{t\rightarrow+\infty}\left(h_\alpha(t)\right)^d{\rm Var}_{\nu_\gamma}\left(\mathbb{E}_{\eta_0}V(\eta_t)\right)
=\frac{\left(\overline{V}^\prime(\gamma)\right)^2{\rm Var}_{\nu_\gamma}\left(\eta(\mathbf{0})\right)f_1^{d, \alpha}(\mathbf{0})}{\left(2\beta^\prime(\gamma)\right)^{\frac{d}{\min\{2, \alpha\}}}},
\]
where $h_\alpha$ is defined as in Lemma \ref{lemma LCLT}.
\end{lemma}

In the simplest case where $V(\eta)=\eta(\mathbf{0})-\gamma$, $c(k)=c_0k$ for all $k\geq 0$ and some $c_0>0$, Lemma \ref{lemma relaxation to equibrium} is a corollary of Lemma \ref{lemma LCLT}. In detail, in this case all particles perform independent long-range random walks with generator $c_0\mathcal{G}^{d, \alpha}$ and $\nu$ is the Poisson distribution with mean $\frac{\beta}{c_0}$. Hence $\beta(\gamma)=c_0\gamma, \beta^\prime(\gamma)=c_0, \overline{V}^\prime(\gamma)=1$ and
\[
\mathbb{E}_{\eta_0}\left(V(\eta_t)\right)=\sum_{y}\eta_0(y)p_{c_0t}^{d, \alpha}(y, \mathbf{0})-\gamma.
\]
Then, according to the strong Markov property of the random walk,
\begin{align*}
{\rm Var}_{\nu_\gamma}\left(\mathbb{E}_{\eta_0}V(\eta_t)\right)&={\rm Var}_{\nu_\gamma}\left(\eta(\mathbf{0})\right)\sum_y\left(p_{c_0t}^{d, \alpha}(y, \mathbf{0})\right)^2\\
&={\rm Var}_{\nu_\gamma}\left(\eta(\mathbf{0})\right)p_{2c_0t}^{d, \alpha}(\mathbf{0}, \mathbf{0}).
\end{align*}
By Lemma \ref{lemma LCLT},
\[
\lim_{t\rightarrow+\infty}\left(h_\alpha(t)\right)^dp_{2c_0t}^{d, \alpha}(\mathbf{0}, \mathbf{0})=f_{2c_0}^{d, \alpha}(\mathbf{0})=f_{1}^{d, \alpha}(\mathbf{0})\frac{1}{\left(2c_0\right)^{\frac{d}{\min\{2, \alpha\}}}}
\]
and hence Lemma \ref{lemma relaxation to equibrium} holds in this special case.

In general cases, Lemma \ref{lemma relaxation to equibrium} follows from an argument similar to that in the nearest-neighbor case given in \cite{Janvresse1999}, where a cut-off technique and a spectral gap estimation play the key roles. In this paper we avoid repeating details similar to those in \cite{Janvresse1999}. Here we only give an intuitive explanation of Lemma \ref{lemma relaxation to equibrium} from the view of the hydrodynamic limit. For any $t\geq 0$ and $N\geq 1$, when the process $\{\eta_s\}_{s\geq 0}$ starts from $\nu_\gamma$, the fluctuation density field $\mathcal{Y}_t^s$ of the process is defined as
\[
\mathcal{Y}_t^s(du)=\frac{1}{\left(h_\alpha(s)\right)^{\frac{d}{2}}}\sum_{x\in \mathbb{Z}^d}\left(\eta_{ts}(x)-\gamma\right)\delta_{\frac{x}{h_\alpha(s)}}(du),
\]
where $\delta_a(du)$ is the Dirac-measure concentrated on $a$. According to an argument similar to that given in \cite{Zhao2025}, as $s\rightarrow+\infty$, $\mathcal{Y}_t^s$ converges weakly, under the Skorohod topology, to $\mathcal{Y}_t$, where $\mathcal{Y}_t$ is the generalized Ornstein-Uhlenbeck process satisfying
\begin{equation}\label{equ O-U}
d\mathcal{Y}_t=\beta^\prime(\gamma)\hat{\mathcal{G}}^{d, \alpha}\mathcal{Y}_tdt+\sqrt{-\beta(\gamma)\hat{\mathcal{G}}^{d, \alpha}}d\mathbb{W}_t,
\end{equation}
where the initial state $\mathcal{Y}_0$ satisfies that $\mathcal{Y}_0(f)$ follows the normal distribution with mean $0$ and variance
\[
{\rm Var}_{\nu_\gamma}\left(\eta(\mathbf{0})\right)\int_{\mathbb{R}^d}f^2(u)du=\frac{\beta(\gamma)}{\beta^\prime(\gamma)}\int_{\mathbb{R}^d}f^2(u)du
\]
for any Schwartz function $f$ on $\mathbb{R}^d$ and $\{\mathbb{W}_t\}_{t\geq 0}$ is the time-space white noise such that $\{\mathbb{W}_t(f)\}_{t\geq 0}$ has the same distribution as that of
$\left\{\sqrt{\int_{\mathbb{R}^d}f^2(u)du}\mathcal{W}_t\right\}_{t\geq 0}$ for any Schwartz function $f$ on $\mathbb{R}^d$. Solving Equation \eqref{equ O-U} by the integral factor method, we have
\[
{\rm Cov}_{\nu_\gamma}\left(\mathcal{Y}_t(f), \mathcal{Y}_s(g)\right)=\frac{\beta(\gamma)}{\beta^\prime(\gamma)}\int_{\mathbb{R}^d}\mathcal{P}_{t-s}^{d, \alpha}f(u)g(u)du
\]
for any $0\leq s<t$ and Schwartz functions $f, g$ on $\mathbb{R}^d$, where $\{\mathcal{P}_{t}^{d, \alpha}\}_{t\geq 0}$ is the semi-group with respect to $\beta^\prime(\gamma)\hat{\mathcal{G}}^{d, \alpha}$, i.e.,
\[
\mathcal{P}_{t}^{d, \alpha}f(u)=\mathbb{E}\left(f\left(\hat{X}^{d, \alpha}_{\beta^\prime(\gamma)t}\right)\Big|\hat{X}_0^{d, \alpha}=u\right).
\]
For sufficiently small $\epsilon>0$, let $h_\epsilon(u)=\frac{1}{(2\epsilon)^d}1_{[-\epsilon, \epsilon]^d}(u)$, then
\begin{align*}
\lim_{\epsilon\rightarrow 0}{\rm Cov}_{\nu_\gamma}\left(\mathcal{Y}_t(h_\epsilon), \mathcal{Y}_0(h_\epsilon)\right)
&=\frac{\beta(\gamma)}{\beta^\prime(\gamma)}\lim_{\epsilon\rightarrow 0}\frac{1}{(2\epsilon)^d}\int_{[-\epsilon, \epsilon]^d}\frac{1}{(2\epsilon)^d}
\mathbb{P}\left(\hat{X}^{d, \alpha}_{\beta^\prime(\gamma)t}\in [-\epsilon, \epsilon]^d\Big|\hat{X}_0^{d, \alpha}=u\right)du\\
&=\frac{\beta(\gamma)}{\beta^\prime(\gamma)}f^{d, \alpha}_{\beta^\prime(\gamma)t}(\mathbf{0})
={\rm Var}_{\nu_\gamma}\left(\eta(\mathbf{0})\right)\frac{f^{d, \alpha}_{1}(\mathbf{0})}{\left(\beta^\prime(\gamma)t\right)^\frac{d}{\min\{2, \alpha\}}}.
\end{align*}
In the view of the hydrodynamic limit of the zero-range process, for large $s$ and small $\epsilon$, we usually can replace $V(\eta_{ts})$ by
\[
\overline{V}\left(\frac{1}{\left(2\epsilon h_\alpha(s)\right)^d}\sum_{x\in \left[-2\epsilon h_\alpha(s), 2\epsilon h_\alpha(s)\right]^d}\eta_{ts}(x)\right)
\]
with small error. Specially, we usually can replace $\eta_{ts}(\mathbf{0})-\gamma$ by
\[
\frac{1}{\left(2\epsilon h_\alpha(s)\right)^d}\sum_{x\in \left[-2\epsilon h_\alpha(s), 2\epsilon h_\alpha(s)\right]^d}(\eta_{ts}(x)-\gamma)
=(h_\alpha(s))^{-\frac{d}{2}}\mathcal{Y}_t^s(h_\epsilon)
\]
with small error. Then, intuitively, for large $s$ and small $\epsilon$,
\begin{align*}
{\rm Cov}_{\nu_\gamma}\left(\eta_0(\mathbf{0}), \eta_{ts}(\mathbf{0})\right)
&\approx (h_\alpha(s))^{-d}{\rm Cov}_{\nu_\gamma}\left(\mathcal{Y}_0^s(h_\epsilon), \mathcal{Y}_t^s(h_\epsilon)\right)\\
&\approx (h_\alpha(s))^{-d}{\rm Cov}_{\nu_\gamma}\left(\mathcal{Y}_0(h_\epsilon), \mathcal{Y}_t(h_\epsilon)\right)\\
&\approx (h_\alpha(s))^{-d}{\rm Var}_{\nu_\gamma}\left(\eta(\mathbf{0})\right)\frac{f^{d, \alpha}_{1}(\mathbf{0})}{\left(\beta^\prime(\gamma)t\right)^\frac{d}{\min\{2, \alpha\}}}.
\end{align*}
Since $\nu_\gamma$ is reversible for the zero-range process, we have
\[
{\rm Var}_{\nu_\gamma}\left(\mathbb{E}_{\eta_0}(\eta_s(\mathbf{0})-\gamma)\right)={\rm Cov}_{\nu_\gamma}\left(\eta_0(\mathbf{0}), \eta_{2s}(\mathbf{0})\right)
\]
and then
\[
{\rm Var}_{\nu_\gamma}\left(\mathbb{E}_{\eta_0}(\eta_s(\mathbf{0})-\gamma)\right)\approx (h_\alpha(s))^{-d}{\rm Var}_{\nu_\gamma}\left(\eta(\mathbf{0})\right)\frac{f^{d, \alpha}_{1}(\mathbf{0})}{\left(2\beta^\prime(\gamma)\right)^\frac{d}{\min\{2, \alpha\}}}
\]
when $s$ is large, which gives an intuitive explanation of Lemma \ref{lemma relaxation to equibrium} when $V(\eta)=\eta(\mathbf{0})-\gamma$.

For general $V$ with polynomial bound, by the Lagrange's mean-value theorem,
\begin{align*}
V(\eta_{tN})& \approx \overline{V}\left(\frac{1}{\left(2\epsilon h_\alpha(s)\right)^d}\sum_{x\in \left[-2\epsilon h_\alpha(s), 2\epsilon h_\alpha(s)\right]^d}\eta_{ts}(x)\right)\\
&\approx \overline{V}(\gamma)+\overline{V}^\prime(\gamma)\left(\frac{1}{\left(2\epsilon h_\alpha(s)\right)^d}\sum_{x\in \left[-2\epsilon h_\alpha(s), 2\epsilon h_\alpha(s)\right]^d}(\eta_{ts}(x)-\gamma)\right)\\
&\approx \overline{V}(\gamma)+\overline{V}^\prime(\gamma)\left(\eta_{ts}(\mathbf{0})-\gamma\right)
\end{align*}
and then
\[
{\rm Var}_{\nu_\gamma}\left(\mathbb{E}_{\eta_0}V(\eta_{ts})\right)\approx \left(\overline{V}^\prime(\gamma)\right)^2{\rm Var}_{\nu_\gamma}\left(\mathbb{E}_{\eta_0}(\eta_{ts}(\mathbf{0})-\gamma)\right),
\]
which gives an intuitive explanation of Lemma \ref{lemma relaxation to equibrium}.

\section{Proof of Theorem \ref{theorem 2.1 main theorem d=1}}\label{section four}
\subsection{The case of $\alpha<1$}\label{subsection 4.1}
In this section, we prove Theorem \ref{theorem 2.1 main theorem d=1}, so throughout this section we assume that $d=1$ and $V: \mathbb{N}^{\mathbb{Z}^1}\rightarrow \mathbb{R}$ with polynomial bound satisfies that $\overline{V}(\gamma)=0$ and $\overline{V}^\prime(\gamma)\neq 0$. We first deal with case where $\alpha<1$, which is relative easy according to the Kipnis-Varadhan theorem.

\proof[Proof of Theorem \ref{theorem 2.1 main theorem d=1} in the case $\alpha<1$]

According to the bilinear property of the covariance and Lemma \ref{lemma relaxation to equibrium}, when $\alpha<1$,
\begin{equation}\label{equ 4.1}
\lim_{N\rightarrow+\infty}{\rm Var}_{\nu_\gamma}\left(\frac{1}{\sqrt{N}}\int_0^{tN}V(\eta_s)ds\right)=t\sigma^2_\gamma(V)<+\infty.
\end{equation}

Since $\nu_\gamma$ is an invariant distribution of $\{\eta_t\}_{t\geq 0}$, $\{\eta_t\}_{t\geq 0}$ starting from $\nu_\gamma$ is stationary. Furthermore, according to an argument similar to that in the proof of Theorem 5.3 of \cite{Komorowski2012}, this stationary process is ergodic. Consequently, Theorem \ref{theorem 2.1 main theorem d=1} in the case $\alpha<1$ follows from \eqref{equ 4.1} and Theorem 1.8 of \cite{Kipnis1986} (Kipnis-Varadhan theorem).
\qed

Note that, the argument in the above proof applies in all cases where $\sigma_\gamma(V)<+\infty$. As a result, we also complete the proofs of Theorem \ref{theorem 2.2 main theorem d=2} in the case $\alpha<2$ and Theorem \ref{theorem 2.3 main theorem d geq 3}.

\subsection{The case of $\alpha=1$}\label{subsection 4.2}
Now we deal with the case of $\alpha=1$. For any $N\geq 1$ and $\eta\in \mathbb{N}^{\mathbb{Z}^1}$, we define
\[
H_N(\eta)=\sum_{x\in \mathbb{Z}^1}\left(\eta(x)-\gamma\right)g_N(x),
\]
where
\[
g_N(x)=\int_0^{+\infty}e^{-\frac{1}{N}s}p_s^{1, 1}(\mathbf{0}, x)ds.
\]
For any $t\geq 0$, we further define
\[
M_t^N=H_N(\eta_t)-H_N(\eta_0)-\int_0^t\mathcal{L}^{1, 1}H_N(\eta_s)ds,
\]
then, according to Dynkin's martingale formula, $\{M_t^N\}_{t\geq 0}$ is a martingale. For any $x\in \mathbb{Z}^1$, according to Kolmogorov-Chapman equation and integral-by-parts formula,
\begin{align*}
\sum_{y\neq x}\left(g_N(y)-g_N(x)\right)\|y-x\|_2^{-(d+\alpha)}
&=\int_0^{+\infty}e^{-s/N}\frac{d}{ds}p_s^{1, 1}(\mathbf{0}, x)ds\\
&=e^{-s/N}p_s^{1, 1}(\mathbf{0}, x)\Big|_0^{+\infty}+\frac{1}{N}\int_0^{+\infty}e^{-s/N}p_s^{1, 1}(\mathbf{0}, x)ds.
\end{align*}
Then
\[
\mathcal{L}^{1, 1}H_N(\eta)=\sum_{x\in \mathbb{Z}^1}\left(c(\eta(x))-\beta(\gamma)\right)\left(-p_0^{1, 1}(\mathbf{0}, x)+\frac{1}{N}g_N(x)\right)
\]
and hence
\[
M_t^N=H_N(\eta_t)-H_N(\eta_0)+\int_0^t\left(c(\eta_s(\mathbf{0}))-\beta(\gamma)\right)ds-Z_t^N,
\]
where
\[
Z_t^N=\frac{1}{N}\int_0^t\sum_{x\in \mathbb{Z}^1}\left(c(\eta_s(x))-\beta(\gamma)\right)g_N(x)ds.
\]
Consequently,
\begin{align}\label{equ 4.2}
&\frac{1}{\sqrt{N\log N}}\int_0^{tN}\left(c(\eta_s(\mathbf{0}))-\beta(\gamma)\right)ds\\
&=\frac{1}{\sqrt{N\log N}}M_{tN}^N-\frac{1}{\sqrt{N\log N}}\left(H_N(\eta_{tN})-H_N(\eta_0)-Z_{tN}^N\right). \notag
\end{align}
We first show that the term $\frac{1}{\sqrt{N\log N}}\left(H_N(\eta_{tN})-H_N(\eta_0)-Z_{tN}^N\right)$ converges weakly to $0$ as $N\rightarrow+\infty$.
\begin{lemma}\label{lemma 4.1}
For any $t\geq 0$, when $\eta_0$ is distributed with $\nu_\gamma$,
\[
\lim_{N\rightarrow+\infty}\frac{1}{\sqrt{N\log N}}\left(H_N(\eta_{tN})-H_N(\eta_0)-Z_{tN}^N\right)=0
\]
in probability.
\end{lemma}

\proof
According to the invariance of $\nu_\gamma$ and Cauchy-Schwarz inequality, we only need to show that
\begin{equation}\label{equ 4.3}
\lim_{N\rightarrow+\infty}{\rm Var}_{\nu_\gamma}\left(\frac{1}{\sqrt{N\log N}}\sum_{x\in \mathbb{Z}^1}\left(\eta_0(x)-\gamma\right)g_N(x)\right)=0
\end{equation}
and
\begin{equation}\label{equ 4.4}
\lim_{N\rightarrow+\infty}{\rm Var}_{\nu_\gamma}\left(\frac{1}{\sqrt{N\log N}}\sum_{x\in \mathbb{Z}^1}\left(c(\eta_0(x))-\beta(\gamma)\right)g_N(x)\right)=0.
\end{equation}
Equations \eqref{equ 4.3} and \eqref{equ 4.4} follow from the same argument and hence here we only check \eqref{equ 4.3}. Since $\{\eta_0(x)\}_{x\in \mathbb{Z}^1}$ are i.i.d. under $\nu_\gamma$,
\begin{align*}
{\rm Var}_{\nu_\gamma}\left(\frac{1}{\sqrt{N\log N}}\sum_{x\in \mathbb{Z}^1}\left(\eta_0(x)-\gamma\right)g_N(x)\right)
&=\frac{1}{N\log N}{\rm Var}_{\nu_\gamma}\left(\eta_0(\mathbf{0})\right)\sum_{x\in \mathbb{Z}^1}\left(g_N(x)\right)^2\\
&=O(1)\frac{1}{N\log N}\int_0^{+\infty}\int_0^{+\infty}e^{-(t+s)/N}p_{s+t}^{1, 1}\left(\mathbf{0}, \mathbf{0}\right)dsdt\\
&=O(1)\frac{1}{N\log N}\int_0^{+\infty}ue^{-u/N}p_u^{1,1}\left(\mathbf{0}, \mathbf{0}\right)du.
\end{align*}
By Lemma \ref{lemma LCLT}, there exists $M<+\infty$ such that
\[
up_u^{1,1}\left(\mathbf{0}, \mathbf{0}\right)\leq 2f_1^{1, 1}(\mathbf{0})
\]
for any $u\geq M$. Therefore,
\begin{align*}
\frac{1}{N\log N}\int_0^{+\infty}ue^{-u/N}p_u^{1,1}\left(\mathbf{0}, \mathbf{0}\right)du&
\leq \frac{O(1)}{N\log N}\left(\int_0^M udu+2f_1^{1, 1}(\mathbf{0})\int_M^{+\infty}e^{-u/N}du\right)\\
&=O\left(\frac{1}{\log N}\right).
\end{align*}
As a result, Equation \eqref{equ 4.3} holds and the proof is complete.
\qed

Now we give the central limit theorem of $\{\frac{1}{\sqrt{N\log N}}M_{tN}^N\}_{t\geq 0}$.
\begin{lemma}\label{lemma 4.2}
Let $\eta_0$ be distributed with $\nu_\gamma$. For any given $T>0$, as $N\rightarrow+\infty$,
\[
\left\{\frac{1}{\sqrt{N\log N}}M_{tN}^N:~0\leq t\leq T\right\}
\]
converges weakly, with respect to the Skorohod topology, to
\[
\left\{\sqrt{2f_1^{1, 1}(\mathbf{0})\beta(\gamma)}\mathcal{W}_t:~0\leq t\leq T\right\}.
\]
\end{lemma}

\proof
We denote by $\left\{\langle M^N \rangle_t\right\}_{t\geq 0}$ the quadratic variation process of $\{M^N_t\}_{t\geq 0}$. By Dynkin's martingale formula,
\begin{equation}\label{equ 4.5}
\langle M^N\rangle_t=\int_0^t\sum_{x\in \mathbb{Z}^1}\sum_{y\neq x}c\left(\eta_s(x)\right)\|y-x\|_2^{-(d+\alpha)}\left(g_N(y)-g_N(x)\right)^2ds.
\end{equation}
According to Theorem 1.4 in Chapter 7 of \cite{Ethier1986}, to complete this proof, we only need to show that
\begin{equation}\label{equ 4.6}
\lim_{N\rightarrow+\infty}\frac{1}{N\log N}\langle M^N \rangle_{tN}=2f_1^{1, 1}(\mathbf{0})\beta(\gamma)t
\end{equation}
in probability for any $t\geq 0$. By Chebyshev's inequality, to prove Equation \eqref{equ 4.6}, we only need to show that
\begin{equation}\label{equ 4.7}
\lim_{N\rightarrow+\infty}\frac{1}{N\log N}\mathbb{E}_{\nu_\gamma}\langle M^N \rangle_{tN}=2f_1^{1, 1}(\mathbf{0})\beta(\gamma)t
\end{equation}
and
\begin{equation}\label{equ 4.8}
\lim_{N\rightarrow+\infty}{\rm Var}_{\nu_\gamma}\left(\frac{1}{N\log N}\langle M^N \rangle_{tN}\right)=0.
\end{equation}
We claim that
\begin{equation}\label{equ 4.9}
\lim_{N\rightarrow+\infty}\frac{1}{\log N}\sum_{x\in \mathbb{Z}^1}\sum_{y\neq x}\|y-x\|_2^{-(d+\alpha)}\left(g_N(y)-g_N(x)\right)^2=2f_1^{1, 1}(\mathbf{0}).
\end{equation}
We check \eqref{equ 4.9} at the end of this proof. By \eqref{equ 1.3}, \eqref{equ 4.5} and \eqref{equ 4.9}, Equation \eqref{equ 4.7} holds. According to the invariance of $\nu_\gamma$ and the bilinear property of the covariance, to prove \eqref{equ 4.8}, we only need to show that
\begin{equation}\label{equ 4.10}
\lim_{t\rightarrow+\infty}\sup_{x,y\in \mathbb{Z}^1}\left|{\rm Cov}_{\nu_\gamma}\left(c\left(\eta_0(x)\right), c\left(\eta_t(y)\right)\right)\right|=0.
\end{equation}
According to the Markov property of $\{\eta_t\}_{t\geq 0}$,
\[
{\rm Cov}_{\nu_\gamma}\left(c\left(\eta_0(x)\right), c\left(\eta_t(y)\right)\right)
=\mathbb{E}_{\nu_\gamma}\Big(\big(c(\eta_0(x))-\beta(\gamma)\big)\mathbb{E}_{\eta_0}\big(c(\eta_t(y))-\beta(\gamma)\big)\Big).
\]
Then, by Lemma \ref{lemma relaxation to equibrium} and Cauchy-Schwarz inequality,
\[
\sup_{x,y\in \mathbb{Z}^1}\left|{\rm Cov}_{\nu_\gamma}\left(c\left(\eta_0(x)\right), c\left(\eta_t(y)\right)\right)\right|=O(t^{-1})
\]
and hence \eqref{equ 4.10} holds.

At last, we only need to check \eqref{equ 4.9}. According to Kolmogorov-Chapman equation,
\begin{align}\label{equ 4.11}
&\sum_{x\in \mathbb{Z}^1}\sum_{y\neq x}\|y-x\|_2^{-(d+\alpha)}\left(g_N(y)-g_N(x)\right)^2 \notag\\
&=-2\sum_{x\in \mathbb{Z}^1}\int_0^{+\infty}p_s^{1, 1}(\mathbf{0}, x)e^{-s/N}\left(\int_0^{+\infty}e^{-t/N}\frac{d}{dt}p_t^{1, 1}(\mathbf{0}, x)dt\right)ds \notag\\
&=-2\sum_{x\in \mathbb{Z}^1}\int_0^{+\infty}p_s^{1, 1}(\mathbf{0}, x)e^{-s/N}\left(-p_0(\mathbf{0}, x)+\frac{1}{N}\int_0^{+\infty}p_t^{1, 1}(\mathbf{0}, x)e^{-t/N}dt\right)ds\\
&=2\int_0^{+\infty}e^{-s/N}p_s^{1,1}(\mathbf{0}, \mathbf{0})ds-\frac{2}{N}\int_0^{+\infty}\int_0^{+\infty}e^{-(t+s)/N}p_{t+s}^{1, 1}(\mathbf{0}, \mathbf{0})dtds. \notag
\end{align}
As we have shown in the proof of Lemma \ref{lemma 4.1},
\[
\lim_{N\rightarrow+\infty}\frac{1}{N\log N}\int_0^{+\infty}\int_0^{+\infty}e^{-(t+s)/N}p_{t+s}^{1, 1}(\mathbf{0}, \mathbf{0})dtds=0.
\]
For any $\epsilon>0$, by Lemma \ref{lemma LCLT}, there exists $M_2=M_2(\epsilon)<+\infty$ such that
\[
(1-\epsilon)f_1^{1, 1}(\mathbf{0})\leq sp^{1,1}_s(\mathbf{0}, \mathbf{0})\leq (1+\epsilon)f_1^{1, 1}(\mathbf{0})
\]
when $s\geq M_2$. Therefore, for sufficiently large $N$,
\begin{align*}
\int_0^{+\infty}e^{-s/N}p_s^{1, 1}(\mathbf{0}, \mathbf{0})ds&\geq (1-\epsilon)f_1^{1, 1}(\mathbf{0})\int_{M_2}^{+\infty}e^{-s/N}\frac{1}{s}ds\\
&\geq f_1^{1, 1}(\mathbf{0})(1-\epsilon)\int_{M_2/N}^\epsilon e^{-\theta}\frac{1}{\theta}d\theta\\
&\geq f_1^{1, 1}(\mathbf{0})(1-\epsilon)e^{-\epsilon}\left(\log N+\log\epsilon-\log M_2\right).
\end{align*}
As a result,
\[
\liminf_{N\rightarrow+\infty}\frac{1}{\log N}\int_0^{+\infty}e^{-s/N}p_s^{1, 1}(\mathbf{0}, \mathbf{0})ds\geq f_1^{1, 1}(\mathbf{0})(1-\epsilon)e^{-\epsilon}.
\]
Let $\epsilon\rightarrow 0$ in the above inequality and then
\begin{equation}\label{equ 4.12}
\liminf_{N\rightarrow+\infty}\frac{1}{\log N}\int_0^{+\infty}e^{-s/N}p_s^{1, 1}(\mathbf{0}, \mathbf{0})ds\geq f_1^{1, 1}(\mathbf{0}).
\end{equation}
For sufficiently large $N$,
\begin{align*}
\int_0^{+\infty}e^{-s/N}p_s^{1, 1}(\mathbf{0}, \mathbf{0})ds&\leq M_2+(1+\epsilon)f_1^{1, 1}(\mathbf{0})\int_{M_2}^{+\infty}e^{-s/N}\frac{1}{s}ds\\
&\leq M_2+(1+\epsilon)f_1^{1, 1}(\mathbf{0})\int^1_{\frac{M_2}{N}}\frac{1}{\theta}d\theta+(1+\epsilon)f_1^{1, 1}(\mathbf{0})\int_1^{+\infty}e^{-\theta}d\theta\\
&=O(1)+(1+\epsilon)f_1^{1, 1}(\mathbf{0})\left(\log N-\log M_2\right)
\end{align*}
and hence
\[
\limsup_{N\rightarrow+\infty}\frac{1}{\log N}\int_0^{+\infty}e^{-s/N}p_s^{1, 1}(\mathbf{0}, \mathbf{0})ds\leq (1+\epsilon)f_1^{1, 1}(\mathbf{0}).
\]
Let $\epsilon\rightarrow 0$ in the above inequality and then
\begin{equation}\label{equ 4.13}
\limsup_{N\rightarrow+\infty}\frac{1}{\log N}\int_0^{+\infty}e^{-s/N}p_s^{1, 1}(\mathbf{0}, \mathbf{0})ds\leq f_1^{1, 1}(\mathbf{0}).
\end{equation}
By \eqref{equ 4.12} and \eqref{equ 4.13}, we have
\begin{equation}\label{equ 4.14}
\lim_{N\rightarrow+\infty}\frac{1}{\log N}\int_0^{+\infty}e^{-s/N}p_s^{1, 1}(\mathbf{0}, \mathbf{0})ds=f_1^{1, 1}(\mathbf{0}).
\end{equation}
Equation \eqref{equ 4.9} follows from \eqref{equ 4.11} and \eqref{equ 4.14}. Since \eqref{equ 4.9} holds, the proof is complete.
\qed

At last, we prove Theorem \ref{theorem 2.1 main theorem d=1} in the case $\alpha=1$.

\proof[Proof of Theorem \ref{theorem 2.1 main theorem d=1} in the case $\alpha=1$]
By Equation \eqref{equ 4.2} and Lemmas \ref{lemma 4.1}, \ref{lemma 4.2},
\[
\frac{1}{\sqrt{N\log N}}\left(\int_0^{t_1N}(c(\eta_s(\mathbf{0}))-\beta(\gamma))ds, \int_0^{t_2N}(c(\eta_s(\mathbf{0}))-\beta(\gamma))ds, \ldots, \int_0^{t_mN}(c(\eta_s(\mathbf{0}))-\beta(\gamma))ds\right)
\]
converges weakly to
\[
\sqrt{2f_1^{1, 1}(\mathbf{0})\beta(\gamma)}\left(\mathcal{W}_{t_1}, \ldots, \mathcal{W}_{t_m}\right)
\]
as $N\rightarrow+\infty$. Let
\[
\Xi(\eta)=\left(\eta(\mathbf{0})-\gamma\right)-\frac{1}{\beta^\prime(\gamma)}\left(c(\eta(0))-\beta(\gamma)\right),
\]
then $\Xi: \mathbb{N}^{\mathbb{Z}^d}\rightarrow \mathbb{R}$ is a local function with polynomial bound and satisfies
\[
\overline{\Xi}(\gamma)=\overline{\Xi}^\prime(\gamma)=0.
\]
Therefore, by Lemma \ref{lemma relaxation to equibrium},
\begin{equation}\label{equ 4.15}
\lim_{s\rightarrow+\infty}s{\rm Var}_{\nu_\gamma}\left(\mathbb{E}_{\eta_0}\Xi(\eta_s)\right)=0.
\end{equation}
According to the invariance of $\nu_\gamma$ and the bilinear property of the covariance,
\[
{\rm Var}_{\nu_\gamma}\left(\int_0^{tN}\Xi(\eta_s)ds\right)
=2\int_0^{tN}\left(\int_0^s\mathbb{E}_{\nu_\gamma}\left(\Xi(\eta_0)\mathbb{E}_{\eta_0}\Xi(\eta_u)\right)du\right)ds.
\]
As we have explained in Section \ref{section three}, since $\nu_\gamma$ is reversible for $\{\eta_t\}_{t\geq 0}$,
\[
\mathbb{E}_{\nu_\gamma}\left(\Xi(\eta_0)\mathbb{E}_{\eta_0}\Xi(\eta_u)\right)={\rm Var}_{\nu_\gamma}\left(\mathbb{E}_{\eta_0}\Xi(\eta_{u/2})\right).
\]
Then, by \eqref{equ 4.15}, for any $\epsilon>0$, there exists $M_3=M_3(\epsilon)<+\infty$ such that
\[
\mathbb{E}_{\nu_\gamma}\left(\Xi(\eta_0)\mathbb{E}_{\eta_0}\Xi(\eta_u)\right)\leq \frac{\epsilon}{u}
\]
for all $u\geq M_3$. Therefore, for sufficiently large $N$,
\begin{align*}
{\rm Var}_{\nu_\gamma}\left(\int_0^{tN}\Xi(\eta_s)ds\right)&\leq 2\int_{M_3}^{tN}\left(\int_{M_3}^s\frac{\epsilon}{u}du\right)ds+O(N)\\
&=\epsilon O(N\log N)
\end{align*}
and hence
\[
\limsup_{N\rightarrow+\infty}\frac{1}{N\log N}{\rm Var}_{\nu_\gamma}\left(\int_0^{tN}\Xi(\eta_s)ds\right)\leq \epsilon.
\]
Let $\epsilon\rightarrow+0$ in the above inequality and then
\begin{equation}\label{equ 4.16}
\lim_{N\rightarrow+\infty}\frac{1}{N\log N}{\rm Var}_{\nu_\gamma}\left(\int_0^{tN}\Xi(\eta_s)ds\right)
=\lim_{N\rightarrow+\infty}\mathbb{E}_{\nu_\gamma}\left(\left(\frac{1}{\sqrt{N\log N}}\int_0^{tN}\Xi(\eta_s)ds\right)^2\right)=0.
\end{equation}
By \eqref{equ 4.16},
\[
\frac{1}{\sqrt{N\log N}}\left(\int_0^{t_1N}(\eta_s(\mathbf{0})-\gamma)ds, \int_0^{t_2N}(\eta_s(\mathbf{0})-\gamma)ds, \ldots, \int_0^{t_mN}(\eta_s(\mathbf{0})-\gamma)ds\right)
\]
converges weakly to
\[
\frac{1}{\beta^\prime(\gamma)}\sqrt{2f_1^{1, 1}(\mathbf{0})\beta(\gamma)}\left(\mathcal{W}_{t_1}, \ldots, \mathcal{W}_{t_m}\right)
\]
as $N\rightarrow+\infty$. Let
\[
\Upsilon(\eta)=V(\eta)-\overline{V}^\prime(\gamma)(\eta(\mathbf{0})-\gamma),
\]
then $\Upsilon: \mathbb{N}^{\mathbb{Z}^d}\rightarrow \mathbb{R}$ is a local function with polynomial bound and satisfies
\[
\overline{\Upsilon}(\gamma)=\overline{\Upsilon}^\prime(\gamma)=0.
\]
According to the same argument as that leading to \eqref{equ 4.16}, we have
\begin{equation}\label{equ 4.17}
\lim_{N\rightarrow+\infty}\mathbb{E}_{\nu_\gamma}\left(\left(\frac{1}{\sqrt{N\log N}}\int_0^{tN}\Upsilon(\eta_s)ds\right)^2\right)=0.
\end{equation}
By \eqref{equ 4.17},
\[
\frac{1}{\sqrt{N\log N}}\left(\int_0^{t_1N}V(\eta_s)ds, \int_0^{t_2N}V(\eta_s)ds, \ldots, \int_0^{t_mN}V(\eta_s)ds\right)
\]
converges weakly to
\[
\frac{\overline{V}^\prime(\gamma)}{\beta^\prime(\gamma)}\sqrt{2f_1^{1, 1}(\mathbf{0})\beta(\gamma)}\left(\mathcal{W}_{t_1}, \ldots, \mathcal{W}_{t_m}\right)
\]
as $N\rightarrow+\infty$ and the proof is complete.
\qed

\subsection{The case of $\alpha>1$}\label{subsection 4.3}
In this subsection, we deal with the case of $\alpha>1$. We first introduce some notations. For any $0\leq s\leq t, \eta\in \mathbb{N}^{\mathbb{Z}^d}$ and $x\in \mathbb{Z}^1$, we define
\[
v(t, x)=\int_0^{\beta^\prime(\gamma)t}p_r^{1, \alpha}(\mathbf{0}, x)dr
\]
and
\[
\Phi_s^t(\eta)=\sum_{x\in \mathbb{Z}^1}\left(\eta(x)-\gamma\right)v(t-s, x).
\]
We further define
\[
\mathcal{M}_s^t=\Phi_s^t(\eta_s)-\Phi_0^t(\eta_0)-\int_0^s(\partial_u+\mathcal{L}^{1, \alpha})\Phi_u^t(\eta_u)du,
\]
then $\{\Phi_s^t\}_{0\leq s\leq t}$ is a martingale according to the Dynkin's martingale formula. According to the definition of $\mathcal{L}^{1, \alpha}$ and Kolmogorov-Chapman equation,
\begin{align*}
\mathcal{L}^{1, \alpha}\Phi_u^t(\eta_u)&=\sum_{x\in \mathbb{Z}^1}\sum_{y\neq x}\|y-x\|_2^{-(d+\alpha)}\left(c(\eta_u(x))-\beta(\gamma)\right)\left(v(t-u, y)-v(t-u, x)\right)\\
&=\sum_{x\in \mathbb{Z}^1}\left(c(\eta_u(x))-\beta(\gamma)\right)\int_0^{\beta^\prime(\gamma)(t-u)}\frac{d}{d\theta}p_\theta^{1, \alpha}(\mathbf{0}, x)d\theta\\
&=\sum_{x\in \mathbb{Z}^1}\left(c(\eta_u(x))-\beta(\gamma)\right)\left(p^{1, \alpha}_{\beta^\prime(\gamma)(t-u)}(\mathbf{0}, x)-p_0^{1, \alpha}(\mathbf{0}, x)\right).
\end{align*}
According to the definition of $\Phi_s^t$,
\[
\partial_u\Phi_u^t(\eta_u)=-\sum_{x\in \mathbb{Z}^1}\left(\eta_u(x)-\gamma\right)p^{1, \alpha}_{\beta^\prime(\gamma)(t-u)}(\mathbf{0}, x)\beta^\prime(\gamma)
\]
and hence
\[
(\partial_u+\mathcal{L}^{1, \alpha})\Phi_u^t(\eta_u)=\sum_{x\in \mathbb{Z}^1}p^{1, \alpha}_{\beta^\prime(\gamma)(t-u)}(\mathbf{0}, x)\zeta(x, \eta_u)-\left(c\left(\eta_u(\mathbf{0})\right)-\beta(\gamma)\right),
\]
where
\[
\zeta(x, \eta)=c(\eta(x))-\beta(\gamma)-\beta^\prime(\gamma)\left(\eta(x)-\gamma\right).
\]
As a result,
\[
\mathcal{M}_s^t=\Phi_s^t(\eta_s)-\Phi_0^t(\eta_0)+\int_0^s\left(c(\eta_u(\mathbf{0}))-\beta(\gamma)\right)du
-\int_0^s\sum_{x\in \mathbb{Z}^1}p^{1, \alpha}_{\beta^\prime(\gamma)(t-u)}(\mathbf{0}, x)\zeta(x, \eta_u)du.
\]
Then, since $\Phi_t^t(\eta)=0$, we have
\begin{equation}\label{equ 4.18 martingale decomposition}
\frac{1}{\Lambda_{1, \alpha}(N)}\int_0^{tN}\left(c(\eta_u(\mathbf{0}))-\beta(\gamma)\right)du=
\frac{1}{\Lambda_{1, \alpha}(N)}\left(\mathcal{M}_{tN}^{tN}+\Phi_0^{tN}(\eta_0)+\mathcal{R}_t^N\right),
\end{equation}
where
\[
\mathcal{R}_t^N=\int_0^{tN}\sum_{x\in \mathbb{Z}^1}p^{1, \alpha}_{\beta^\prime(\gamma)(tN-u)}(\mathbf{0}, x)\zeta(x, \eta_u)du.
\]
We first show that the remainder term $\frac{1}{\Lambda_{1, \alpha}(N)}\mathcal{R}_t^N$ converges weakly to $0$ as $N\rightarrow+\infty$.

\begin{lemma}\label{lemma 4.4}
Let $\alpha>1$. For any $t\geq 0$, when $\eta_0$ is distributed with $\nu_\gamma$,
\[
\lim_{N\rightarrow+\infty}\frac{1}{\Lambda_{1, \alpha}(N)}\mathcal{R}_t^N=0
\]
in $L^2$.
\end{lemma}

\proof
We only deal with the case of $1<\alpha<2$ since other two cases follow from the same argument. According to the invariance of $\nu_\gamma$,
\begin{align*}
&\mathbb{E}_{\nu_\gamma}\left(\left(\mathcal{R}_t^N\right)^2\right)\\
&=2\sum_{x\in \mathbb{Z}^1}\sum_{y\in \mathbb{Z}^1}\int_0^{tN}p^{1, \alpha}_{\beta^\prime(\gamma)(tN-u)}(\mathbf{0}, x)
\left(\int_0^up^{1, \alpha}_{\beta^\prime(\gamma)(tN+\theta-u)}(\mathbf{0}, y)
\mathbb{E}_{\nu_\gamma}\left(\zeta(x, \eta_0)\mathbb{E}_{\eta_0}\zeta(y, \eta_\theta)\right)d\theta\right)du.
\end{align*}
Since $\nu_\gamma$ is reversible for $\{\eta_t\}_{t\geq 0}$, we have
\[
\mathbb{E}_{\nu_\gamma}\left(\zeta(x, \eta_0)\mathbb{E}_{\eta_0}\zeta(y, \eta_\theta)\right)
=\mathbb{E}_{\nu_\gamma}\left(\mathbb{E}_{\eta_0}\zeta(x, \eta_{\theta/2})\mathbb{E}_{\eta_0}\zeta(y, \eta_{\theta/2})\right).
\]
Since $\overline{\zeta(x)}(\gamma)=\overline{\zeta(x)}^\prime(\gamma)=0$, for any $\epsilon>0$, there exists $M_5=M_5(\epsilon)<+\infty$ such that
\[
\theta^{-1/\alpha}\sup_{x, y\in \mathbb{Z}^1}\mathbb{E}_{\nu_\gamma}\left(\zeta(x, \eta_{\theta/2})\mathbb{E}_{\eta_0}\zeta(y, \eta_{\theta/2})\right)
\leq \epsilon
\]
when $\theta\geq M_5$ by Lemma \ref{lemma relaxation to equibrium} and Cauchy-Schwarz inequality. Consequently, for sufficiently large $N$ and $1<\alpha<2$,
\begin{align*}
\mathbb{E}_{\nu_\gamma}\left(\left(\mathcal{R}_t^N\right)^2\right)
&\leq O(N)+O(1)\epsilon\int_{M_5}^{tN}\left(\int_{M_5}^u\theta^{-1/\alpha}d\theta\right)du\\
&=O(1)\epsilon N^{2-\frac{1}{\alpha}}=O(1)\epsilon \Lambda_{1, \alpha}^2(N).
\end{align*}
Therefore,
\[
\limsup_{N\rightarrow+\infty}\mathbb{E}_{\nu_\gamma}\left(\left(\frac{1}{\Lambda_{1, \alpha}(N)}\mathcal{R}_t^N\right)^2\right)\leq \epsilon.
\]
Let $\epsilon\rightarrow 0$ in the above inequality and then the proof is complete.
\qed

The following lemma gives the CLT of the term $\frac{1}{\Lambda_{1, \alpha}(N)}\Phi_0^{tN}(\eta_0)$ in \eqref{equ 4.18 martingale decomposition}.

\begin{lemma}\label{lemma 4.5}
Let $\alpha>1$ and $\eta_0$ be distributed with $\nu_\gamma$. For any $0\leq t_1<t_2<\ldots<t_m$, as $N\rightarrow+\infty$,
\[
\frac{1}{\Lambda_{1, \alpha}(N)}\left(\Phi_0^{t_1N}(\eta_0), \ldots, \Phi_0^{t_mN}(\eta_0)\right)
\]
converges weakly to $\left(\Psi_1, \ldots, \Psi_m\right)$, which is a $\mathbb{R}^m$-valued Gaussian random variable such that $\mathbb{E}\Psi_j=0$ and
\begin{align*}
&{\rm Cov}\left(\Psi_i, \Psi_j\right)=\\
&\begin{cases}
\frac{\alpha^2}{(\alpha-1)(2\alpha-1)}f_1^{1, \alpha}(\mathbf{0})\beta(\gamma)\left(\beta^\prime(\gamma)\right)^{1-\frac{1}{\alpha}}
\left((t_i+t_j)^{2-\frac{1}{\alpha}}-t_i^{2-\frac{1}{\alpha}}-t_j^{2-\frac{1}{\alpha}}\right) & \text{~if~}1<\alpha<2,\\
\frac{4}{3}f_1^{1, \alpha}(\mathbf{0})\beta(\gamma)\left(\beta^\prime(\gamma)\right)^{\frac{1}{2}}
\left((t_i+t_j)^{\frac{3}{2}}-t_i^{\frac{3}{2}}-t_j^{\frac{3}{2}}\right) & \text{~if~}\alpha\geq 2
\end{cases}
\end{align*}
for any $1\leq i\leq j\leq m$.
\end{lemma}

\proof
We only deal with the case of $1<\alpha<2$ since the other two cases follow from similar arguments. We denote by $H^N$ the characteristic function of
\[
\frac{1}{\Lambda_{1, \alpha}(N)}\left(\Phi_0^{t_1N}(\eta_0), \ldots, \Phi_0^{t_mN}(\eta_0)\right),
\]
i.e.,
\[
H^N(\vec{u})=\mathbb{E}_{\nu_\gamma}\exp\left\{\frac{\sqrt{-1}}{\Lambda_{1, \alpha}(N)}\sum_{k=1}^mu_k\Phi_0^{t_kN}(\eta_0)\right\}
\]
for any $\vec{u}=(u_1, \ldots, u_m)\in \mathbb{R}^m$. To complete this proof, we only need to show that
\begin{equation}\label{equ 4.19}
\lim_{N\rightarrow+\infty}H_N(\vec{u})=e^{-\frac{1}{2}\vec{u}\Sigma \vec{u}^T},
\end{equation}
where $\Sigma_{i, j}={\rm Cov}(\Psi_i, \Psi_j)$ for all $1\leq i, j\leq m$.

Since $\{\eta(x)\}_{x\in \mathbb{Z}^1}$ are i.i.d. under $\nu_\gamma$, we have
\[
H_N(\vec{u})=\exp\left\{\sum_{x\in \mathbb{Z}^1}\mathcal{C}(x)\right\},
\]
where
\[
\mathcal{C}(x)=\log \left(\mathbb{E}_{\nu_\gamma}\exp\left\{\left(\eta_0(x)-\gamma\right)\frac{\sqrt{-1}}{\Lambda_{1, \alpha}(N)}
\left(\sum_{k=1}^mu_kv(t_kN, x)\right)\right\}\right).
\]
Since $\mathbb{E}_{\nu_\gamma}\left(\eta_0(x)-\gamma\right)=0$, according to the fact that $\log(1+x)=x+o(x), e^x=1+x+\frac{x^2}{2}+o(x^2)$ and \eqref{equ 1.4}, to prove \eqref{equ 4.19} we only need to check that
\begin{align}\label{equ 4.20}
&\lim_{N\rightarrow+\infty}\frac{1}{\Lambda_{1, \alpha}^2(N)}\sum_{x\in \mathbb{Z}^1}\int_0^{aN}p^{1, \alpha}_\theta(\mathbf{0}, x)d\theta\int_0^{bN}p^{1, \alpha}_\theta(\mathbf{0}, x)d\theta\\
&=\frac{\alpha^2}{(\alpha-1)(2\alpha-1)}f_1^{1, \alpha}(\mathbf{0})\left((a+b)^{2-\frac{1}{\alpha}}-a^{2-\frac{1}{\alpha}}-b^{2-\frac{1}{\alpha}}\right) \notag
\end{align}
for any $a, b>0$. According to the Markov property of the long-range random walk,
\begin{align*}
\sum_{x\in \mathbb{Z}^1}\int_0^{aN}p^{1, \alpha}_\theta(\mathbf{0}, x)d\theta\int_0^{bN}p^{1, \alpha}_\theta(\mathbf{0}, x)d\theta
&=\int_0^{aN}\left(\int_0^{bN}p^{1, \alpha}_{s+\theta}(\mathbf{0}, \mathbf{0})ds\right)d\theta.
\end{align*}
When $1<\alpha<2$, by Lemma \ref{lemma LCLT}, for any $\epsilon>0$ there exists $M_6=M_6(\epsilon)<+\infty$ such that
\[
(1-\epsilon)f_1^{1, \alpha}(\mathbf{0})\leq r^{1/\alpha}p^{1, \alpha}_r(\mathbf{0}, \mathbf{0})\leq (1+\epsilon)f_1^{1, \alpha}(\mathbf{0})
\]
for $r\geq M_6$. Therefore, for sufficiently large $N$,
\begin{align*}
&\int_0^{aN}\left(\int_0^{bN}p^{1, \alpha}_{s+\theta}(\mathbf{0}, \mathbf{0})ds\right)d\theta\\
&\geq (1-\epsilon)f_1^{1, \alpha}(\mathbf{0})\int_{M_6}^{aN}\left(\int_0^{bN}(s+\theta)^{-1/\alpha}ds\right)d\theta\\
&=\frac{\alpha(1-\epsilon)f_1^{1, \alpha}(\mathbf{0})}{\alpha-1}\int_{M_6}^{aN}\left(\left(bN+\theta\right)^{1-1/\alpha}-\theta^{1-1/\alpha}\right)d\theta\\
&=\frac{\alpha^2(1-\epsilon)f_1^{1, \alpha}(\mathbf{0})N^{2-1/\alpha}}{(\alpha-1)(2-\alpha)}
\left((a+b)^{2-1/\alpha}-(b+M_6/N)^{2-1/\alpha}-a^{2-1/\alpha}+(M_6/N)^{2-1/\alpha}\right)
\end{align*}
and hence
\begin{align*}
&\liminf_{N\rightarrow+\infty}\frac{1}{\Lambda_{1, \alpha}^2(N)}\sum_{x\in \mathbb{Z}^1}\int_0^{aN}p^{1, \alpha}_\theta(\mathbf{0}, x)d\theta\int_0^{bN}p^{1, \alpha}_\theta(\mathbf{0}, x)d\theta\\
&\geq \frac{\alpha^2(1-\epsilon)}{(\alpha-1)(2\alpha-1)}f_1^{1, \alpha}(\mathbf{0})\left((a+b)^{2-\frac{1}{\alpha}}-a^{2-\frac{1}{\alpha}}-b^{2-\frac{1}{\alpha}}\right).
\end{align*}
Let $\epsilon\rightarrow0$ in the above inequality and then
\begin{align}\label{equ 4.21}
&\liminf_{N\rightarrow+\infty}\frac{1}{\Lambda_{1, \alpha}^2(N)}\sum_{x\in \mathbb{Z}^1}\int_0^{aN}p^{1, \alpha}_\theta(\mathbf{0}, x)d\theta\int_0^{bN}p^{1, \alpha}_\theta(\mathbf{0}, x)d\theta\notag\\
&\geq \frac{\alpha^2}{(\alpha-1)(2\alpha-1)}f_1^{1, \alpha}(\mathbf{0})\left((a+b)^{2-\frac{1}{\alpha}}-a^{2-\frac{1}{\alpha}}-b^{2-\frac{1}{\alpha}}\right).
\end{align}
For sufficiently large $N$,
\begin{align*}
&\int_0^{aN}\left(\int_0^{bN}p^{1, \alpha}_{s+\theta}(\mathbf{0}, \mathbf{0})ds\right)d\theta\\
&\leq O(N)+(1+\epsilon)f_1^{1, \alpha}(\mathbf{0})\int_{M_6}^{aN}\left(\int_0^{bN}(s+\theta)^{-1/\alpha}ds\right)d\theta\\
&=O(N)+\frac{\alpha(1+\epsilon)f_1^{1, \alpha}(\mathbf{0})}{\alpha-1}\int_{M_6}^{aN}\left(\left(bN+\theta\right)^{1-1/\alpha}-\theta^{1-1/\alpha}\right)d\theta\\
&=O(N)+\frac{\alpha^2(1+\epsilon)f_1^{1, \alpha}(\mathbf{0})N^{2-1/\alpha}}{(\alpha-1)(2-\alpha)}
\left((a+b)^{2-1/\alpha}-(b+M_6/N)^{2-1/\alpha}-a^{2-1/\alpha}+(M_6/N)^{2-1/\alpha}\right)
\end{align*}
and hence
\begin{align*}
&\limsup_{N\rightarrow+\infty}\frac{1}{\Lambda_{1, \alpha}^2(N)}\sum_{x\in \mathbb{Z}^1}\int_0^{aN}p^{1, \alpha}_\theta(\mathbf{0}, x)d\theta\int_0^{bN}p^{1, \alpha}_\theta(\mathbf{0}, x)d\theta\\
&\leq \frac{\alpha^2(1+\epsilon)}{(\alpha-1)(2\alpha-1)}f_1^{1, \alpha}(\mathbf{0})\left((a+b)^{2-\frac{1}{\alpha}}-a^{2-\frac{1}{\alpha}}-b^{2-\frac{1}{\alpha}}\right).
\end{align*}
Let $\epsilon\rightarrow0$ in the above inequality and then
\begin{align}\label{equ 4.22}
&\limsup_{N\rightarrow+\infty}\frac{1}{\Lambda_{1, \alpha}^2(N)}\sum_{x\in \mathbb{Z}^1}\int_0^{aN}p^{1, \alpha}_\theta(\mathbf{0}, x)d\theta\int_0^{bN}p^{1, \alpha}_\theta(\mathbf{0}, x)d\theta\notag\\
&\leq \frac{\alpha^2}{(\alpha-1)(2\alpha-1)}f_1^{1, \alpha}(\mathbf{0})\left((a+b)^{2-\frac{1}{\alpha}}-a^{2-\frac{1}{\alpha}}-b^{2-\frac{1}{\alpha}}\right).
\end{align}
Equation \eqref{equ 4.20} follows from \eqref{equ 4.21} and \eqref{equ 4.22}. Hence, the proof is complete.
\qed

Now we deal with the CLT of the term $\frac{1}{\Lambda_{1, \alpha}(N)}\mathcal{M}_{tN}^{tN}$  in \eqref{equ 4.18 martingale decomposition}. According to the definition of $\{\eta_t\}_{t\geq 0}$, there exist independent Poisson processes $\left\{N_{x, y}(t):~t\geq 0\right\}_{x, y\in \mathbb{Z}^1}$  with rate $1$ such that
\[
\eta_t(x)=\eta_0(x)+\sum_{y\neq x}N_{y, x}\left(\int_0^tc(\eta_s(y))\|y-x\|_2^{-(1+\alpha)}ds\right)-\sum_{y\neq x}N_{x, y}\left(\int_0^tc(\eta_s(x))\|y-x\|_2^{-(1+\alpha)}ds\right)
\]
for all $x\in \mathbb{Z}^1$ and $t\geq 0$. Let
\[
J_{x, y}(t)=N_{x, y}\left(\int_0^tc(\eta_s(x))\|y-x\|_2^{-(1+\alpha)}ds\right)-\int_0^tc(\eta_s(x))\|y-x\|_2^{-(1+\alpha)}ds,
\]
then $\{J_{x, y}(t)\}_{t\geq 0}$ ia a martingale. Furthermore, for any $x, y, z, w\in \mathbb{Z}^1$, the cross-variation process $\left\{\langle J_{x, y}, J_{z, w} \rangle_t\right\}_{t\geq 0}$ of $\{J_{x, y}(t)\}_{t\geq 0}$ and $\{J_{z, w}(t)\}_{t\geq 0}$ is given by
\[
\langle J_{x, y}, J_{z, w} \rangle_t=
\begin{cases}
c(\eta_t(x))\|y-x\|_2^{-(1+\alpha)} & \text{~if~} (x, y)=(z, w),\\
0 & \text{~otherwise}.
\end{cases}
\]
According to the definition of $\mathcal{M}_s^t$, by It\^{o}'s formula, we have
\begin{equation}\label{equ 4.23}
\mathcal{M}_s^t=\sum_{x\in \mathbb{Z}^1}\sum_{y\neq x}\int_0^sv(t-\theta, y)-v(t-\theta, x)dJ_{x, y}(\theta)
\end{equation}
for any $0\leq s\leq t$. For given $T>0$ and any $N\geq 1$, $\alpha>1$, we denote by $\mathcal{U}^N_\alpha$ the random signed measure on $[0, T]\times \mathbb{R}$ such that
\begin{align*}
\mathcal{U}_\alpha^N(H)=\varpi(\alpha, N)\sum_{x\in \mathbb{Z}^1}\sum_{y\in \mathbb{Z}^1\setminus\{x\}}\int_0^T\left(\iint\limits_{\mathcal{A}(x, y, \alpha, N)}H(s, v)-H(s, u)dudv\right)dJ_{x, y}(Ns)
\end{align*}
for any $H\in C_c^\infty([0, T]\times \mathbb{R})$, where
\[
\mathcal{A}(x, y, \alpha, N)=(x, y)/{h_\alpha(N)}
+\left(-\frac{1}{2h_\alpha(N)}, \frac{1}{2h_\alpha(N)}\right]^2
\]
and
\[
\varpi(\alpha, N)=
\begin{cases}
N^{\frac{3}{2\alpha}} & \text{~if~} 1<\alpha<2,\\
N^{\frac{3}{4}}\left(\log N\right)^{\frac{3}{4}} & \text{~if~} \alpha=2,\\
N^{\frac{3}{4}} & \text{~if~}\alpha>2.
\end{cases}.
\]
For any $u\in \mathbb{R}$, we denote by $u_{\alpha, N}$ the element in
\[
\mathbb{Z}^1/{h_\alpha(N)}:=\left\{\frac{x}{h_\alpha(N)}:~x\in \mathbb{Z}^1\right\}
\]
such that
\[
u-u_{\alpha, N}\in \left(-\frac{1}{2h_\alpha(N)}, \frac{1}{2h_\alpha(N)}\right].
\]
Then, for any $N\geq 1$, $0\leq t, \theta\leq T$ and $u\in \mathbb{R}$, we define
\begin{align*}
b_t^N(\theta, u)&=\left(\int_0^{\beta^\prime(\gamma)(t-\theta)}h_\alpha(N)p_{Ns}^{1, \alpha}(\mathbf{0}, h_{\alpha}(N)u_{\alpha, N})ds\right)1_{\{\theta\leq t\}}\\
&=1_{\{\theta\leq t\}}h_\alpha(N)N^{-1}v\left((t-\theta)N, h_\alpha(N)u_{\alpha, N}\right).
\end{align*}
Then, by \eqref{equ 4.23} and the definition of $\mathcal{U}_\alpha^N$, for any $t\leq T$,
\begin{equation}\label{equ 4.24}
\frac{1}{\Lambda_{1, \alpha}(N)}\mathcal{M}_{tN}^{tN}=\mathcal{U}_\alpha^N(b_t^N).
\end{equation}
We further define $\mathcal{U}_\alpha$ as the Gaussian time-space white noise on $[0, T]\times \mathbb{R}$ such that $\mathcal{U}_\alpha(H)$ follows the normal distribution with mean $0$ and variance $\varrho(\alpha, H)$ for any $H\in C_c^{+\infty}\left([0, T]\times \mathbb{R}\right)$, where
\[
\varrho(\alpha, H)=\beta(\gamma)\int_0^T\left(\iint_{\mathbb{R}^2}\|v-u\|_2^{-(1+\alpha)}\left(H(s, v)-H(s, u)\right)^2dudv\right)ds
\]
when $1<\alpha<2$,
\[
\varrho(\alpha, H)=\beta(\gamma)\int_0^T\left(\int_{\mathbb{R}}\left(\partial_uH(s, u)\right)^2du\right)ds
\]
when $\alpha=2$ and
\[
\varrho(\alpha, H)=2\beta(\gamma)\int_0^T\left(\int_{\mathbb{R}}\left(\partial_uH(s, u)\right)^2du\right)ds\left(\sum_{k=1}^{+\infty}\frac{1}{k^{\alpha-1}}\right)
\]
when $\alpha>2$.

The following lemma gives the CLT of $\{\mathcal{U}_\alpha^N\}_{N\geq 1}$.
\begin{lemma}\label{lemma 4.6}
Let $\alpha>1$ and $\eta_0$ be distributed with $\nu_\gamma$. For any $H\in C_c^\infty([0, T]\times \mathbb{R})$, as $N\rightarrow+\infty$,
$\mathcal{U}_\alpha^N(H)$ converges weakly to $\mathcal{U}_\alpha(H)$.
\end{lemma}

We give a detailed proof of Lemma \ref{lemma 4.6} in the case of $1<\alpha<2$. For the other two cases, we give an outline of the proofs.

\proof[Proof of Lemma \ref{lemma 4.6} in the case of $1<\alpha<2$]
Assume that $1<\alpha<2$. For any $h\in C_c^\infty(\mathbb{R})$ and $0\leq s, t\leq T, u\in \mathbb{R}$, we define
\[
\widetilde{h}_t(s, u)=h(u)1_{\{s\leq t\}}.
\]
Then, to complete the proof, we only need to show that $\left\{\mathcal{U}^N_\alpha(\widetilde{h}_t)\right\}_{0\leq t\leq T}$ converges weakly to
\[
\left\{\sqrt{\beta(\gamma)\iint_{\mathbb{R}^2}\|v-u\|_2^{-(1+\alpha)}\left(h(v)-h(u)\right)^2dudv}\mathcal{W}_t\right\}_{0\leq t\leq T}
\]
as $N\rightarrow+\infty$. We denote by $\{\langle \mathcal{U}^N_\alpha(\widetilde{h}_\cdot) \rangle_t\}_{0\leq t\leq T}$ the quadratic variation process of $\left\{\mathcal{U}^N_\alpha(\widetilde{h}_t)\right\}_{0\leq t\leq T}$. To complete the proof, we only need to show that
\begin{equation}\label{equ 4.25}
\lim_{N\rightarrow+\infty}\langle \mathcal{U}^N_\alpha(\widetilde{h}_\cdot) \rangle_t=\left(\beta(\gamma)\iint_{\mathbb{R}^2}\|v-u\|_2^{-(1+\alpha)}\left(h(v)-h(u)\right)^2dudv\right)t
\end{equation}
in probability. To check \eqref{equ 4.25}, we only need to show that
\begin{equation}\label{equ 4.26}
\lim_{N\rightarrow+\infty}\mathbb{E}_{\nu_\gamma}\langle \mathcal{U}^N_\alpha(\widetilde{h}_\cdot) \rangle_t
=\left(\beta(\gamma)\iint_{\mathbb{R}^2}\|v-u\|_2^{-(1+\alpha)}\left(h(v)-h(u)\right)^2dudv\right)t
\end{equation}
and
\begin{equation}\label{equ 4.27}
\lim_{N\rightarrow+\infty}{\rm Var}_{\nu_\gamma}\left(\langle \mathcal{U}^N_\alpha(\widetilde{h}_\cdot) \rangle_t\right)=0.
\end{equation}
According to the definition of $\mathcal{U}_\alpha^N$, we have
\begin{align*}
\langle \mathcal{U}^N_\alpha(\widetilde{h}_\cdot) \rangle_t
&=N^{\frac{3}{\alpha}}N\sum_{x}\sum_{y\neq x}\int_0^tc(\eta_{sN}(x))\|y-x\|_2^{-(1+\alpha)}\left(\iint\limits_{\mathcal{A}(x, y, \alpha, N)}h(v)-h(u)dudv\right)^2ds\\
&=(1+o(1))N^{1+\frac{3}{\alpha}}\sum_{x}\sum_{y\neq x}\int_0^tc(\eta_{sN}(x))\|y-x\|_2^{-(1+\alpha)}\left(h(y/{N^{\frac{1}{\alpha}}})-h(x/{N^{\frac{1}{\alpha}}})\right)^2N^{-\frac{4}{\alpha}}ds\\
&=\frac{1+o(1)}{N^{\frac{1}{2\alpha}}}\sum_{x}\sum_{y\neq x}\int_0^tc(\eta_{sN}(x))\left\|\frac{y}{N^{\frac{1}{\alpha}}}-\frac{x}{N^{\frac{1}{\alpha}}}\right\|_2^{-(1+\alpha)}
\left(h(y/{N^{\frac{1}{\alpha}}})-h(x/{N^{\frac{1}{\alpha}}})\right)^2ds.
\end{align*}
Consequently, \eqref{equ 4.26} follows from \eqref{equ 1.3} and the definition of Riemann Integral.

According to the bilinear property of the covariance, to prove \eqref{equ 4.27}, we only need to check a `$1<\alpha<2$'-version of \eqref{equ 4.10}. According to Lemma \ref{lemma relaxation to equibrium} and the same argument as that given in the proof of Lemma \ref{lemma 4.2}, when $1<\alpha<2$,
\[
\sup_{x,y\in \mathbb{Z}^1}\left|{\rm Cov}_{\nu_\gamma}\left(c\left(\eta_0(x)\right), c\left(\eta_t(y)\right)\right)\right|=O(t^{-1/\alpha})
\]
and hence \eqref{equ 4.27} holds. In conclusion, the proof is complete.
\qed

\proof[Proof of Lemma \ref{lemma 4.6} in the case $\alpha\geq 2$]
According to the same argument as that given in the proof of the case $1<\alpha<2$, we only need to show that
\begin{equation}\label{equ 4.28}
\lim_{N\rightarrow+\infty}\mathbb{E}_{\nu_\gamma}\langle \mathcal{U}^N_\alpha(\widetilde{h}_\cdot) \rangle_t
=\left(\beta(\gamma)\int_{\mathbb{R}}\left(\partial_uh(u)\right)^2du
\right)t
\end{equation}
when $\alpha=2$,
\begin{equation}\label{equ 4.29}
\lim_{N\rightarrow+\infty}\mathbb{E}_{\nu_\gamma}\langle \mathcal{U}^N_\alpha(\widetilde{h}_\cdot) \rangle_t
=\left(2\beta(\gamma)\int_{\mathbb{R}}\left(\partial_uh(u)\right)^2du\left(\sum_{k=1}^{+\infty}\frac{1}{k^{\alpha-1}}\right)\right)t
\end{equation}
when $\alpha>2$ and then check a `$\alpha\geq 2$'-version of \eqref{equ 4.27}.  When $\alpha>2$, for sufficiently small $\epsilon$ and $|y-x|\leq \sqrt{N}\epsilon$,
\[
h(y/\sqrt{N})-h(x/\sqrt{N})=\left(\partial_uh(x/\sqrt{N})+o(1)\right)|y-x|/\sqrt{N}
\]
and hence, for sufficiently large $N$,
\begin{align*}
&N^{\frac{3}{2}}N\sum_{x}\sum_{|y-x|\leq \sqrt{N}\epsilon}\int_0^t\|y-x\|_2^{-(1+\alpha)}\left(\iint\limits_{\mathcal{A}(x, y, \alpha, N)}h(v)-h(u)dudv\right)^2ds\\
&=(1+o(1))N^{5/2}\frac{1}{N^2}\sum_{x}t\left(\partial_uh(x/\sqrt{N})\right)^2\sum_{|y-x|\leq \sqrt{N}\epsilon}|y-x|^{-(1+\alpha)}|(y-x)/\sqrt{N}|^2\\
&=2(1+o(1))tN^{-\frac{1}{2}}\left(\sum_{x}\left(\partial_uh(x/\sqrt{N})\right)^2\right)\sum_{k=1}^{\epsilon\sqrt{N}}\frac{1}{k^{\alpha-1}}\\
&=2(1+o(1))\left(\int_{\mathbb{R}}\left(\partial_uh(u)\right)^2du\left(\sum_{k=1}^{+\infty}\frac{1}{k^{\alpha-1}}\right)\right)t.
\end{align*}
Similarly, by the integrability of $h^2$, for sufficiently large $N$,
\begin{align*}
&N^{\frac{3}{2}}N\sum_{x}\sum_{|y-x|\geq \sqrt{N}\epsilon}\int_0^t\|y-x\|_2^{-(1+\alpha)}\left(\iint\limits_{\mathcal{A}(x, y, \alpha, N)}h(v)-h(u)dudv\right)^2ds\\
&=O(1)N^{5/2}\frac{1}{N^2}\sqrt{N}\sum_{k\geq \epsilon\sqrt{N}}\frac{1}{k^{1+\alpha}}=O(1)N^{(2-\alpha)/2}=o(1).
\end{align*}
Then, \eqref{equ 4.29} follows from \eqref{equ 1.3}. When $\alpha=2$, similarly,
\[
h(y/\sqrt{N\log N})-h(x/\sqrt{N\log N})=\left(\partial_uh(x/\sqrt{N\log N})+o(1)\right)|y-x|/\sqrt{N\log N}
\]
when $|y-x|\leq \epsilon\sqrt{N\log N}$. Then, according to the fact that
\[
\sum_{k=1}^{\epsilon\sqrt{N\log N}}\frac{k^2}{k^3}=\frac{1}{2}(1+o(1))\log N,
\]
we have
\begin{align*}
&N^{\frac{3}{2}}(\log N)^{\frac{3}{2}}N\sum_{x}\sum_{|y-x|\leq \sqrt{N\log N}\epsilon}\int_0^t\|y-x\|_2^{-3}\left(\iint\limits_{\mathcal{A}(x, y, \alpha, N)}h(v)-h(u)dudv\right)^2ds\\
&=\frac{2((1+o(1)))N^{5/2}(\log N)^{3/2}}{N^3(\log N)^3}\sum_{x}t\left(\partial_uh(x/\sqrt{N\log N})\right)^2\left(\sum_{k=1}^{\epsilon\sqrt{N\log N}}\frac{k^2}{k^3}\right)\\
&=(2t)(1+o(1))\left(\frac{1}{\sqrt{N\log N}}\sum_{x}t\left(\partial_uh(x/\sqrt{N\log N})\right)^2\right)\left(\frac{1}{\log N}\sum_{k=1}^{\epsilon\sqrt{N\log N}}\frac{k^2}{k^3}\right)\\
&=t(1+o(1))\int_{\mathbb{R}}\left(\partial_uh(u)\right)^2du.
\end{align*}
According to the integrability of $h^2$ and the fact that
\[
\sum_{k=\epsilon\sqrt{N\log N}}^{+\infty}\frac{1}{k^3}=\frac{O(1)}{N\log N},
\]
it is easy to check that
\begin{align*}
&N^{\frac{3}{2}}(\log N)^{\frac{3}{2}}N\sum_{x}\sum_{|y-x|\geq \sqrt{N\log N}\epsilon}\int_0^t\|y-x\|_2^{-3}\left(\iint\limits_{\mathcal{A}(x, y, \alpha, N)}h(v)-h(u)dudv\right)^2ds\\
&=O(1)\frac{1}{\log N}=o(1).
\end{align*}
Then, \eqref{equ 4.28} follows from \eqref{equ 1.3}. To prove a `$\alpha\geq 2$'-version of \eqref{equ 4.27}, by the bilinear property of the covariance, we only need to check a `$\alpha\geq 2$'-version of \eqref{equ 4.10}. According to Lemma \ref{lemma relaxation to equibrium} and the same argument as that given in the proof of Lemma \ref{lemma 4.2}, we have
\[
\sup_{x,y\in \mathbb{Z}^1}\left|{\rm Cov}_{\nu_\gamma}\left(c\left(\eta_0(x)\right), c\left(\eta_t(y)\right)\right)\right|=O(t^{-1/2})
\]
when $\alpha\geq 2$ and hence the `$\alpha\geq 2$'-version of \eqref{equ 4.10} holds. In conclusion, the proof is complete.
\qed

The following lemma gives the CLT of $\frac{1}{\Lambda_{1, \alpha}(N)}\mathcal{M}_{tN}^{tN}$.
\begin{lemma}\label{lemma 4.7}
Let $\alpha>1$ and $\eta_0$ be distributed with $\nu_\gamma$. For any $0\leq t_1<t_2<\ldots<t_m$,
\[
\frac{1}{\Lambda_{1, \alpha}(N)}\left(\mathcal{M}_{t_1N}^{t_1N}, \mathcal{M}_{t_2N}^{t_2N}, \ldots, \mathcal{M}_{t_mN}^{t_mN}\right)
\]
converges weakly, as $N\rightarrow+\infty$, to $\left(\Gamma_1, \ldots, \Gamma_m\right)$, which is a $\mathbb{R}^m$-valued Gaussian random variable such that
$\mathbb{E}\Gamma_{i}=0$ and
\begin{align*}
&{\rm Cov}\left(\Gamma_i, \Gamma_j\right)=\\
&\begin{cases}
\frac{\alpha^2f_1^{1, \alpha}(\mathbf{0})\beta(\gamma)\left(\beta^\prime(\gamma)\right)^{1-\frac{1}{\alpha}}}{(\alpha-1)(2\alpha-1)}
\left(2t_i^{2-\frac{1}{\alpha}}+2t_j^{2-\frac{1}{\alpha}}-(t_j-t_i)^{2-\frac{1}{\alpha}}-(t_i+t_j)^{2-\frac{1}{\alpha}}\right)
& \text{~if~}1<\alpha<2,\\
\frac{4}{3}f_1^{1, \alpha}(\mathbf{0})\beta(\gamma)\left(\beta^\prime(\gamma)\right)^{\frac{1}{2}}\left(2t_i^{\frac{3}{2}}+2t_j^{\frac{3}{2}}
-(t_i-t_j)^{\frac{3}{2}}-(t_i+t_j)^{\frac{3}{2}}\right) & \text{~if~}\alpha\geq 2
\end{cases}
\end{align*}
for all $1\leq i\leq j\leq m$.
\end{lemma}

\proof
We only deal with the case of $1<\alpha<2$ since the other two cases follow from similar arguments. For given $T>0$ and any $t\in [0, T]$, let $b_t^N$ be defined as in \eqref{equ 4.24}. By Lemma \ref{lemma LCLT}, $b_t^N$ converges to $b_t$ in $L^2([0, T]\times \mathbb{R})$ as $N\rightarrow+\infty$, where
\[
b_t(s, u)=\left(\int_0^{\beta^\prime(\gamma)(t-s)}f_\theta^{1, \alpha}(u)d\theta\right)1_{\{s\leq t\}}
\]
for any $(s, u)\in [0, T]\times \mathbb{R}$. We choose $T$ sufficiently large such that $T>t_m$, then by \eqref{equ 4.24} and Lemma \ref{lemma 4.6},
\[
\frac{1}{\Lambda_{1, \alpha}(N)}\left(\mathcal{M}_{t_1N}^{t_1N}, \mathcal{M}_{t_2N}^{t_2N}, \ldots, \mathcal{M}_{t_mN}^{t_mN}\right)
\]
converges weakly to
\[
\left(\mathcal{U}_\alpha(b_{t_1}),\ldots, \mathcal{U}_\alpha(b_{t_m})\right)
\]
as $N\rightarrow+\infty$. Hence, to complete this proof, we only need to show that
\begin{align*}
&{\rm Cov}\left(\mathcal{U}_\alpha(b_{t_i}), \mathcal{U}_\alpha(b_{t_j})\right)\\
&=\frac{\alpha^2}{(\alpha-1)(2\alpha-1)}f_1^{1, \alpha}(\mathbf{0})\beta(\gamma)\left(\beta^\prime(\gamma)\right)^{1-\frac{1}{\alpha}}
\left(2t_i^{2-\frac{1}{\alpha}}+2t_j^{2-\frac{1}{\alpha}}-(t_j-t_i)^{2-\frac{1}{\alpha}}-(t_i+t_j)^{2-\frac{1}{\alpha}}\right)
\end{align*}
for any $i\leq j$. According to the definition of $b_t$ and Kolmogrov-Chapman equation, for any $i\leq j$ and $s\leq t_i$,
\begin{align*}
&\iint_{\mathbb{R}^2}|v-u|^{-1-\alpha}\left(b_{t_i}(s, v)-b_{t_i}(s,u)\right)\left(b_{t_j}(s, v)-b_{t_j}(s, u)\right)dudv\\
&=(-2)\lim_{\epsilon\rightarrow 0}\int_{\mathbb{R}}b_{t_i}(s, u)\left(\int_\epsilon^{(t_j-s)\beta^\prime(\gamma)}\left(\int_\mathbb{R}|v-u|^{-1-\alpha}\left(f_\theta^{1, \alpha}(v)-f_\theta^{1, \alpha}(u)\right)dv\right)d\theta\right)du\\
&=(-2)\lim_{\epsilon\rightarrow 0}\int_{\mathbb{R}}b_{t_i}(s, u)\left(\int_\epsilon^{(t_j-s)\beta^\prime(\gamma)}\left(\frac{d}{d\theta}f_\theta^{1, \alpha}(u)\right)d\theta\right)du\\
&=(-2)\lim_{\epsilon\rightarrow 0}\int_{\mathbb{R}}b_{t_i}(s, u)\left(f^{1, \alpha}_{(t_j-s)\beta^\prime(\gamma)}(u)-f_\epsilon^{1, \alpha}(u)\right)du\\
&=(-2)\lim_{\epsilon\rightarrow 0}\int_0^{(t_i-s)\beta^\prime(\gamma)}\left(\int_{\mathbb{R}}f_\theta^{1, \alpha}(u)f_{(t_j-s)\beta^\prime(\gamma)}^{1, \alpha}(u)du-\int_{\mathbb{R}}f_\theta^{1, \alpha}(u)f_\epsilon^{1, \alpha}(u)du\right)d\theta\\
&=(-2)\lim_{\epsilon\rightarrow 0}\int_0^{(t_i-s)\beta^\prime(\gamma)}\left(f_{\theta+(t_j-s)\beta^\prime(\gamma)}^{1, \alpha}(\mathbf{0})-f_{\theta+\epsilon}^{1, \alpha}(\mathbf{0})\right)d\theta\\
&=2\int_0^{(t_i-s)\beta^\prime(\gamma)}f_{\theta}^{1, \alpha}(\mathbf{0})-f_{\theta+(t_j-s)\beta^\prime(\gamma)}^{1, \alpha}(\mathbf{0})d\theta\\
&=2f_{1}^{1, \alpha}(\mathbf{0})\int_0^{(t_i-s)\beta^\prime(\gamma)}\theta^{-\frac{1}{\alpha}}
-(\theta+(t_j-s)\beta^\prime(\gamma))^{-\frac{1}{\alpha}}d\theta.
\end{align*}
Consequently, for $i\leq j$,
\begin{align*}
&{\rm Cov}\left(\mathcal{U}_\alpha(b_{t_i}), \mathcal{U}_\alpha(b_{t_j})\right)\\
&=2\beta(\gamma)f_{1}^{1, \alpha}(\mathbf{0})\int_0^{t_i}\left(\int_0^{(t_i-s)\beta^\prime(\gamma)}\theta^{-\frac{1}{\alpha}}
-(\theta+(t_j-s)\beta^\prime(\gamma))^{-\frac{1}{\alpha}}d\theta\right)ds\\
&=\frac{\alpha^2}{(\alpha-1)(2\alpha-1)}f_1^{1, \alpha}(\mathbf{0})\beta(\gamma)\left(\beta^\prime(\gamma)\right)^{1-\frac{1}{\alpha}}
\left(2t_i^{2-\frac{1}{\alpha}}+2t_j^{2-\frac{1}{\alpha}}-(t_j-t_i)^{2-\frac{1}{\alpha}}-(t_i+t_j)^{2-\frac{1}{\alpha}}\right)
\end{align*}
and hence the proof is complete.
\qed

The following lemma gives the CLT of $\frac{1}{\Lambda_{1, \alpha}(N)}\left(\mathcal{M}_{tN}^{tN}, \Phi_0^{tN}(\eta_0)\right)$.

\begin{lemma}\label{lemma 4.8}
Let $\alpha>1$ and $\eta_0$ be distributed with $\nu_\gamma$. For any $0<t_1<t_2<\ldots<t_m$, as $N\rightarrow+\infty$, the joint distribution of
\[
\frac{1}{\Lambda_{1, \alpha}(N)}\left(\Phi_0^{t_1N}(\eta_0), \ldots, \Phi_0^{t_mN}(\eta_0)\right)
\text{~and~}\frac{1}{\Lambda_{1, \alpha}(N)}\left(\mathcal{M}_{t_1N}^{t_1N}, \mathcal{M}_{t_2N}^{t_2N}, \ldots, \mathcal{M}_{t_mN}^{t_mN}\right)
\]
converges weakly to the independent coupling of
\[
\left(\Psi_1, \ldots, \Psi_m\right) \text{~and~} \left(\Gamma_1, \ldots, \Gamma_m\right),
\]
which are defined as in Lemmas \ref{lemma 4.5} and \ref{lemma 4.7} respectively.
\end{lemma}

\proof
We only deal with the case of $1<\alpha<2$ since the other two cases follow from the same argument. By \eqref{equ 4.25}, there exists a subsequence $\{N_k\}_{k\geq 1}$ of $\{N\}_{N\geq 1}$ such that
\[
\lim_{k\rightarrow+\infty}\langle \mathcal{U}^{N_k}_\alpha(\widetilde{h}_\cdot) \rangle_t=\left(\beta(\gamma)\iint_{\mathbb{R}^2}\|v-u\|_2^{-(1+\alpha)}\left(h(v)-h(u)\right)^2dudv\right)t
\]
almost surely for any $t\in \mathbb{Q}$ and $h$ in a given countable dense subset of $C_c^{+\infty}(\mathbb{R})$. Therefore, there exists $\mathcal{K}\subseteq \mathbb{N}^{\mathbb{Z}^1}$ such that $\nu_\gamma(\mathcal{K})=1$ and, under $\mathbb{P}_\eta$ for each $\eta\in \mathcal{K}$, $\langle \mathcal{U}^{N_k}_\alpha(\widetilde{h}_\cdot) \rangle_t$ converges to
\[
\left(\beta(\gamma)\iint_{\mathbb{R}^2}\|v-u\|_2^{-(1+\alpha)}\left(h(v)-h(u)\right)^2dudv\right)t
\]
almost surely for all $t\in \mathbb{Q}$ and $h$ in the aforesaid countable dense subset of $C_c^{+\infty}(\mathbb{R})$ as $k\rightarrow+\infty$. Consequently, under $\mathbb{P}_\eta$ for each $\eta\in \mathcal{K}$,
\[
\frac{1}{\Lambda_{1, \alpha}(N_k)}\left(\mathcal{M}_{t_1N_k}^{t_1N_k}, \mathcal{M}_{t_2N_k}^{t_2N_k}, \ldots, \mathcal{M}_{t_mN_k}^{t_mN_k}\right)
\]
converges weakly to $\left(\Gamma_1, \ldots, \Gamma_m\right)$ as $k\rightarrow+\infty$. Then, by the iterated expectation law,
\[
\frac{1}{\Lambda_{1, \alpha}(N_k)}\left(\Phi_0^{t_1N_k}(\eta_0), \ldots, \Phi_0^{t_mN_k}(\eta_0), \mathcal{M}_{t_1N_k}^{t_1N_k}, \mathcal{M}_{t_2N_k}^{t_2N_k}, \ldots, \mathcal{M}_{t_mN_k}^{t_mN_k}\right)
\]
converges weakly to the independent coupling of
\[
\left(\Psi_1, \ldots, \Psi_m\right) \text{~and~} \left(\Gamma_1, \ldots, \Gamma_m\right)
\]
as $k\rightarrow+\infty$. According to the above argument, any subsequence of
\[
\left\{\frac{1}{\Lambda_{1, \alpha}(N)}\left(\Phi_0^{t_1N}(\eta_0), \ldots, \Phi_0^{t_mN}(\eta_0), \mathcal{M}_{t_1N}^{t_1N}, \mathcal{M}_{t_2N}^{t_2N}, \ldots, \mathcal{M}_{t_mN}^{t_mN}\right)\right\}_{N\geq 1}
\]
have a subsequence that converges weakly to the independent coupling of
\[
\left(\Psi_1, \ldots, \Psi_m\right) \text{~and~} \left(\Gamma_1, \ldots, \Gamma_m\right)
\]
as $k\rightarrow+\infty$. Hence, the proof is complete.
\qed

The following lemma gives the CLT of $\frac{1}{\Lambda_{1, \alpha}(N)}\int_0^{tN}\left(c(\eta_u(\mathbf{0}))-\beta(\gamma)\right)du$.

\begin{lemma}\label{lemma 4.9}
Let $\alpha>1$ and $\eta_0$ be distributed with $\nu_\gamma$. For any $0\leq t_1<t_2<\ldots<t_m$,
\[
\frac{1}{\Lambda_{1, \alpha}(N)}\left(\int_0^{t_1N}\left(c(\eta_u(\mathbf{0}))-\beta(\gamma)\right)du, \int_0^{t_2N}\left(c(\eta_u(\mathbf{0}))-\beta(\gamma)\right)du, \ldots,\int_0^{t_mN}\left(c(\eta_u(\mathbf{0}))-\beta(\gamma)\right)du\right)
\]
converges weakly to
\[
\begin{cases}
\sqrt{\frac{2\alpha^2f_1^{1, \alpha}(\mathbf{0})\beta(\gamma)\left(\beta^\prime(\gamma)\right)^{1-\frac{1}{\alpha}}}{(\alpha-1)(2\alpha-1)}}\left(B_{t_1}^{1-\frac{1}{2\alpha}}, B_{t_2}^{1-\frac{1}{2\alpha}}, \ldots, B_{t_m}^{1-\frac{1}{2\alpha}}\right) & \text{~if~}1<\alpha<2,\\
\sqrt{\frac{8}{3}f_1^{1, \alpha}(\mathbf{0})\beta(\gamma)\left(\beta^\prime(\gamma)\right)^{\frac{1}{2}}}
\left(B_{t_1}^{\frac{3}{4}}, \ldots, B_{t_m}^{\frac{3}{4}}\right) & \text{~if~}\alpha\geq 2
\end{cases}
\]
as $N\rightarrow+\infty$.
\end{lemma}

\proof
We only deal with the case of $1<\alpha<2$, since the other two cases follow from the same argument. By Lemmas \ref{lemma 4.4}, \ref{lemma 4.8} and \eqref{equ 4.18 martingale decomposition}, as $N\rightarrow+\infty$,
\[
\frac{1}{\Lambda_{1, \alpha}(N)}\left(\int_0^{t_1N}\left(c(\eta_u(\mathbf{0}))-\beta(\gamma)\right)du, \int_0^{t_2N}\left(c(\eta_u(\mathbf{0}))-\beta(\gamma)\right)du, \ldots,\int_0^{t_mN}\left(c(\eta_u(\mathbf{0}))-\beta(\gamma)\right)du\right)
\]
converges weakly to the $\mathbb{R}^m$-valued Gaussian random variable
\[
\left(\Psi_1+\Gamma_1, \ldots, \Psi_m+\Gamma_m\right),
\]
where $\left(\Psi_1, \ldots, \Psi_m\right)$ and $\left(\Gamma_1,\ldots, \Gamma_m\right)$ are defined as in Lemmas \ref{lemma 4.5} and \ref{lemma 4.7} respectively and are independently coupled in the same probability space.

According to the definitions of $\Psi_i$ and $\Gamma_i$, when $1<\alpha<2$, we have
\[
\mathbb{E}\left(\Psi_i+\Gamma_i\right)=0
\]
and
\begin{align*}
&{\rm Cov}\left(\Psi_i+\Gamma_i, \Psi_j+\Gamma_j\right)\\
&=\frac{\alpha^2}{(\alpha-1)(2\alpha-1)}f_1^{1, \alpha}(\mathbf{0})\beta(\gamma)\left(\beta^\prime(\gamma)\right)^{1-\frac{1}{\alpha}}
\Big((t_i+t_j)^{2-\frac{1}{\alpha}}-t_i^{2-\frac{1}{\alpha}}-t_j^{2-\frac{1}{\alpha}}
\\
&\text{\quad\quad\quad\quad\quad}+2t_i^{2-\frac{1}{\alpha}}+2t_j^{2-\frac{1}{\alpha}}-(t_j-t_i)^{2-\frac{1}{\alpha}}-(t_i+t_j)^{2-\frac{1}{\alpha}}\Big)\\
&=\frac{\alpha^2}{(\alpha-1)(2\alpha-1)}f_1^{1, \alpha}(\mathbf{0})\beta(\gamma)\left(\beta^\prime(\gamma)\right)^{1-\frac{1}{\alpha}}
\left(t_i^{2-\frac{1}{\alpha}}+t_j^{2-\frac{1}{\alpha}}-(t_j-t_i)^{2-\frac{1}{\alpha}}\right)
\end{align*}
for all $1\leq i\leq j\leq m$. As a result,
\[
\left(\Psi_1+\Gamma_1, \ldots, \Psi_m+\Gamma_m\right)
\]
and
\[
\sqrt{\frac{2\alpha^2f_1^{1, \alpha}(\mathbf{0})\beta(\gamma)\left(\beta^\prime(\gamma)\right)^{1-\frac{1}{\alpha}}}{(\alpha-1)(2\alpha-1)}}\left(B_{t_1}^{1-\frac{1}{2\alpha}}, B_{t_2}^{1-\frac{1}{2\alpha}}, \ldots, B_{t_m}^{1-\frac{1}{2\alpha}}\right)
\]
have the same probability distribution and the proof is complete.
\qed

For any local $V$ with polynomial bound satisfying $\overline{V}(\gamma)=0$, according to the same argument as that given in the proof of Theorem \ref{theorem 2.1 main theorem d=1} in the case $\alpha=1$, we have
\[
\frac{1}{\Lambda_{1, \alpha}(N)}\int_0^{tN}V(\eta_s)ds
=\frac{\overline{V}^\prime(\gamma)}{\Lambda_{1, \alpha}(N)}\int_0^{tN}\left(\eta_s(\mathbf{0})-\gamma\right)ds+\varepsilon_N^1
\]
and
\[
\frac{1}{\Lambda_{1, \alpha}(N)}\int_0^{tN}\left(\eta_s(\mathbf{0})-\gamma\right)ds
=\frac{1}{\Lambda_{1, \alpha}(N)\beta^\prime(\gamma)}\int_0^{tN}\left(c(\eta_s(\mathbf{0}))-\beta(\gamma)\right)ds+\varepsilon_N^2,
\]
where
\[
\lim_{N\rightarrow+\infty}\varepsilon_N^1=\lim_{N\rightarrow+\infty}\varepsilon_N^2=0
\]
in probability. Consequently, by Lemma \ref{lemma 4.9}, we have the following lemma.
\begin{lemma}\label{lemma 4.10}
Let $\eta_0$ be distributed with $\nu_\gamma$ and local $V: \mathbb{N}^{\mathbb{Z}^1}\rightarrow \mathbb{R}$ with polynomial bound satisfy that $\overline{V}(\gamma)=0, \overline{V}^\prime(\gamma)\neq 0$. For any $0<t_1<t_2<\ldots<t_m$,
\[
\frac{1}{\Lambda_{1, \alpha}(N)}\left(\int_0^{t_1N}V(\eta_s)ds, \int_0^{t_2N}V(\eta_s)ds, \ldots, \int_0^{t_mN}V(\eta_s)ds\right)
\]
converges weakly to
\[
\begin{cases}
\sqrt{\frac{2\alpha^2}{(\alpha-1)(2\alpha-1)}\frac{\left(\overline{V}^\prime(\gamma)\right)^2\beta(\gamma)}{\left(\beta^\prime(\gamma)\right)^{1+\frac{1}{\alpha}}}f_1^{1, \alpha}(\mathbf{0})}\left(B_{t_1}^{1-\frac{1}{2\alpha}}, \ldots, B_{t_m}^{1-\frac{1}{2\alpha}}\right) & \text{~if~}1<\alpha<2,\\
\sqrt{\frac{8}{3}\frac{\left(\overline{V}^\prime(\gamma)\right)^2\beta(\gamma)}{\left(\beta^\prime(\gamma)\right)^{\frac{3}{2}}}f_1^{1, \alpha}(\mathbf{0})}\left(B_{t_1}^{\frac{3}{4}}, \ldots, B_{t_m}^{\frac{3}{4}}\right) & \text{~if~}\alpha\geq 2
\end{cases}
\]
as $N\rightarrow+\infty$.
\end{lemma}

At last, we prove Theorem \ref{theorem 2.1 main theorem d=1} in the case $\alpha>1$.

\proof[Proof of Theorem \ref{theorem 2.1 main theorem d=1} in the case $\alpha>1$]
By Lemma \ref{lemma 4.10}, we only need to show that
\[
\left\{\frac{1}{\Lambda_{1, \alpha}(N)}\int_0^{tN}V(\eta_s)ds:~0\leq t\leq T\right\}_{N\geq 1}
\]
are tight under the uniform topology. To check this tightness, according to Corollary 14.9 of \cite{Kallenberg1997} and the invariance of $\nu_\gamma$, we only need to show that there exist $a, b, c>0$ independent of $N, t$ such that
\begin{equation}\label{equ 4.30}
\mathbb{E}_{\nu_\gamma}\left(\left(\frac{1}{\Lambda_{1, \alpha}(N)}\int_0^{tN}V(\eta_s)ds\right)^a\right)
\leq ct^{1+b}
\end{equation}
for all $t\in [0, 1]$ and $N\geq 1$. According to an argument similar to that given in the proof of Theorem \ref{theorem 2.1 main theorem d=1} in the case $\alpha=1$, we have
\[
\mathbb{E}_{\nu_\gamma}\left(\left(\int_0^{tN}V(\eta_s)ds\right)^2\right)
=2\int_0^{tN}\left(\int_0^s{\rm Var}_{\nu_\gamma}\left(\mathbb{E}_{\eta_0}V(\eta_{u/2})\right)du\right)ds.
\]
When $1<\alpha<2$, by Lemma \ref{lemma relaxation to equibrium}, there exists $M_7<+\infty$ such that
\[
\sup_{t\geq 0}t^{1/\alpha}{\rm Var}_{\nu_\gamma}\left(\mathbb{E}_{\eta_0}V(\eta_{t/2})\right)\leq M_7.
\]
Hence, when $1<\alpha<2$,
\begin{align*}
\mathbb{E}_{\nu_\gamma}\left(\left(\int_0^{tN}V(\eta_s)ds\right)^2\right)
&\leq 2M_7\int_0^{tN}\left(\int_0^su^{-1/\alpha}du\right)ds\\
&=\frac{2\alpha^2M_7}{(\alpha-1)(2\alpha-1)}(tN)^{2-\frac{1}{\alpha}}
=\frac{2\alpha^2M_7t^{2-\frac{1}{\alpha}}}{(\alpha-1)(2\alpha-1)}\Lambda^2_{1, \alpha}(N)
\end{align*}
for any $t\geq 0$. Hence, when $1<\alpha<2$, \eqref{equ 4.30} holds with $a=2, b=1-\frac{1}{\alpha}$ and $c=\frac{2\alpha^2M_7}{(\alpha-1)(2\alpha-1)}$.

When $\alpha>2$,  by Lemma \ref{lemma relaxation to equibrium}, there exists $M_8<+\infty$ such that
\[
\sup_{t\geq 0}t^{1/2}{\rm Var}_{\nu_\gamma}\left(\mathbb{E}_{\eta_0}V(\eta_{t/2})\right)\leq M_8.
\]
Then, according to a similar argument, \eqref{equ 4.30} holds with $a=2, b=\frac{1}{2}$ and $c=8M_8/3$.

When $\alpha=2$, by Lemma \ref{lemma relaxation to equibrium}, there exists $M_9<+\infty$ and $1000<M_{10}<+\infty$ such that
\[
\sup_{t\geq M_{10}}\sqrt{t\log t}{\rm Var}_{\nu_\gamma}\left(\mathbb{E}_{\eta_0}V(\eta_{t/2})\right)\leq M_9.
\]
Then for $N\geq 4$ and $0\leq t\leq 1$, when $tN\leq M_{10}$,
\begin{align*}
\mathbb{E}_{\nu_\gamma}\left(\left(\int_0^{tN}V(\eta_s)ds\right)^2\right)
\leq t^2N^2M_{11},
\end{align*}
where
\[
M_{11}=\sup_{0\leq u\leq M_{10}/2}{\rm Var}_{\nu_\gamma}\left(\mathbb{E}_{\eta_0}V(\eta_{u/2})\right)
\leq {\rm Var}_{\nu_\gamma}(V)<+\infty.
\]
Hence, when $N\geq 4, t\in [0, 1]$ and $tN\leq M_{10}$,
\begin{align}\label{equ 4.31}
\mathbb{E}_{\nu_\gamma}\left(\left(\frac{1}{\Lambda_{1, 2}(N)}\int_0^{tN}V(\eta_s)ds\right)^2\right)
& \leq M_{11}t^2N^{1/2}(\log N)^{1/2}\notag\\
&\leq M_{11}M_{12}t^2N^{0.6}\leq M_{10}^{0.6}M_{11}M_{12}t^{1.4},
\end{align}
where
\[
M_{12}=\sup_{N\geq 4}\frac{\sqrt{N}(\log N)^{1/2}}{N^{0.6}}<+\infty.
\]
When $N\geq 4$ and $tN\geq M_{10}$,
\begin{align*}
\mathbb{E}_{\nu_\gamma}\left(\left(\int_0^{tN}V(\eta_s)ds\right)^2\right)
&\leq M_{11}M_{10}^2+2tNM_{10}M_{11}+2M_9\int_{M_{10}}^{tN}\left(\int_{M_{10}}^s\frac{1}{\sqrt{u\log u}}du\right)ds.
\end{align*}
According to the formula of integral-by-parts,
\begin{align*}
\int_{M_{10}}^s\frac{1}{\sqrt{u\log u}}du
&\leq 2s^{1/2}(\log s)^{-1/2}+\int_{M_{10}}^su^{-1/2}(\log u)^{-3/2}du\\
&\leq 2s^{1/2}(\log s)^{-1/2}+ \frac{1}{\log M_{10}}\int_{M_{10}}^su^{-1/2}(\log u)^{-1/2}du
\end{align*}
and hence
\[
\int_{M_{10}}^s\frac{1}{\sqrt{u\log u}}du\leq 2s^{1/2}(\log s)^{-1/2}M_{13},
\]
where $M_{13}=\left(1-\frac{1}{\log M_{10}}\right)^{-1}$. According to a similar argument,
\[
\int_{M_{10}}^{tN}s^{1/2}(\log s)^{-1/2}ds
\leq \frac{2}{3}(tN)^{3/2}(\log (tN))^{-1/2}M_{14},
\]
where $M_{14}=\left(1-\frac{1}{3\log M_{10}}\right)^{-1}$. Therefore, when $N\geq 4, t\in [0, 1]$ and $tN\geq M_{10}$,
\[
2M_9\int_{M_{10}}^{tN}\left(\int_{M_{10}}^s\frac{1}{\sqrt{u\log u}}du\right)ds
\leq M_{15}(tN)^{3/2}(\log (tN))^{-1/2},
\]
where
\[
M_{15}=\frac{8}{3}M_9M_{13}M_{14}.
\]
Then,  when $N\geq 4, t\in [0, 1]$ and $tN\geq M_{10}$,
\[
\mathbb{E}_{\nu_\gamma}\left(\left(\int_0^{tN}V(\eta_s)ds\right)^2\right)
\leq  M_{11}M_{10}^2+2tNM_{10}M_{11}+M_{15}(tN)^{3/2}(\log (tN))^{-1/2}
\]
and hence
\begin{align}\label{equ 4.32}
&\mathbb{E}_{\nu_\gamma}\left(\left(\frac{1}{\Lambda_{1, 2}(N)}\int_0^{tN}V(\eta_s)ds\right)^2\right) \notag\\
&\leq \frac{(\log N)^{1/2}}{N^{3/2}}M_{11}M_{10}^2+\frac{(\log N)^{1/2}}{N^{1/2}}2tM_{10}M_{11}+M_{15}t^{3/2}\left(\frac{\log N}{\log(tN)}\right)^{1/2}\notag\\
&\leq \frac{(\log N)^{1/2}}{N^{0.1}}M_{11}M_{10}^{0.6}t^{1.4}+\frac{2(\log N)^{1/2}}{N^{0.1}}M_{11}M_{10}^{0.6}t^{1.4}
+M_{15}t^{3/2}\left(1+\frac{\log(1/t)}{\log M_{10}}\right)^{1/2}\notag\\
&\leq M_{16}t^{1.4},
\end{align}
where
\[
M_{16}=3M_{10}M_{11}\max_{N\geq 4}\frac{(\log N)^{1/2}}{N^{0.1}}+M_{15}\sup_{t\in[0, 1]}\left(t^{0.1}\left(1+\frac{\log(1/t)}{\log M_{10}}\right)^{1/2}\right).
\]
By \eqref{equ 4.31} and \eqref{equ 4.32}, when $\alpha=2$, \eqref{equ 4.30} holds with $a=2, b=0.4$ and
\[
c=\max\left\{M_{16}, M_{10}^{0.6}M_{11}M_{12}\right\}.
\]
Since \eqref{equ 4.30} holds in all three cases, the proof is complete.
\qed

\begin{remark}\label{remark 4.10}
By Lemma \ref{lemma relaxation to equibrium} and an argument similar to that given in the above proof, when $\alpha=1$ we have the following analogue of \eqref{equ 4.30},
\[
\mathbb{E}_{\nu_\gamma}\left(\left(\frac{1}{\sqrt{N\log N}}\int_0^{tN}V(\eta_s)ds\right)^2\right)
\leq ct^{0.9}
\]
for some $c<+\infty$ and all $t\in [0, 1], N\geq 4$. However, to obtain the tightness of
\[
\left\{\frac{1}{\sqrt{N\log N}}\int_0^{tN}V(\eta_s)ds:~0\leq t\leq T\right\}_{N\geq 1}
\]
under the uniform topology, we need to strengthen $t^{0.9}$ in the above inequality to $t^b$ for some $b>1$ with $\left(\frac{1}{\sqrt{N\log N}}\int_0^{tN}V(\eta_s)ds\right)^2$ being replaced by $\left(\frac{1}{\sqrt{N\log N}}\int_0^{tN}V(\eta_s)ds\right)^a$ for some proper $a>0$. We will work on this problem as a further investigation.
\end{remark}

\section{Proof of Theorem \ref{theorem 2.2 main theorem d=2} in the case $\alpha\geq 2$}\label{section five}

As we have explained in Section \ref{section four}, Theorem \ref{theorem 2.2 main theorem d=2} in the case $\alpha<2$ and Theorem \ref{theorem 2.3 main theorem d geq 3} follow from the same argument as that in the proof of Theorem \ref{theorem 2.1 main theorem d=1} in the case $\alpha<1$. So we only need to prove Theorem \ref{theorem 2.2 main theorem d=2} in the case $\alpha\geq 2$ in this section. Hence, throughout this section we assume that $d=2$ and $\alpha\geq 2$. We only give an outline of the proof of the case of $d=2, \alpha\geq 2$ since all arguments are similar to those given in Subsection \ref{subsection 4.2}.

Let $H_N, g_N, M_t^N, Z_t^N$ be defined as in Subsection \ref{subsection 4.2} except that all $1, 1$ in these notations are replaced by $2, \alpha$. Then, we have
\begin{equation}\label{equ 5.1}
\frac{1}{\Lambda_{2, \alpha}(N)}\int_0^{tN}\left(c(\eta_s(\mathbf{0}))-\beta(\gamma)\right)ds
=\frac{1}{\Lambda_{2, \alpha}(N)}M_{tN}^N-\frac{1}{\Lambda_{2, \alpha}(N)}\left(H_N(\eta_{tN})-H_N(\eta_0)-Z_{tN}^N\right),
\end{equation}
which is an analogue of \eqref{equ 4.2}. Then, we have the following analogue of Lemma \ref{lemma 4.1}.

\begin{lemma}\label{lemma 5.1}
For any $t\geq 0$, when $\eta_0$ is distributed with $\nu_\gamma$,
\[
\lim_{N\rightarrow+\infty}\frac{1}{\Lambda_{2, \alpha}(N)}\left(H_N(\eta_{tN})-H_N(\eta_0)-Z_{tN}^N\right)=0
\]
in probability.
\end{lemma}
\proof
According to an argument similar to that given in the proof of Lemma \ref{lemma 4.1}, we only need to show that
\begin{equation}\label{equ 5.2}
\lim_{N\rightarrow+\infty}\frac{1}{\Lambda^2_{2, \alpha}(N)}\int_0^{+\infty}ue^{-u/N}p_u^{2, \alpha}(\mathbf{0}, \mathbf{0})du=0.
\end{equation}
By Lemma \ref{lemma LCLT},
\[
up_u^{2, \alpha}(\mathbf{0}, \mathbf{0})=
\begin{cases}
O(\frac{1}{\log u}) & \text{~if~} \alpha=2,\\
O(1) & \text{~if~} \alpha>2
\end{cases}
\]
for sufficiently large $u$ and then
\[
\int_0^{+\infty}ue^{-u/N}p_u^{2, \alpha}(\mathbf{0}, \mathbf{0})du=O(N)
\]
for all $\alpha\geq 2$. As a result,
\[
\frac{1}{\Lambda^2_{2, \alpha}(N)}\int_0^{+\infty}ue^{-u/N}p_u^{2, \alpha}(\mathbf{0}, \mathbf{0})du
=\begin{cases}
O\left(\frac{1}{\log(\log N)}\right) & \text{~if~} \alpha=2,\\
O\left(\frac{1}{\log N}\right) & \text{~if~} \alpha>2
\end{cases}
\]
and hence \eqref{equ 5.2} holds.
\qed

The following lemma is an analogue of Lemma \ref{lemma 4.2}.
\begin{lemma}\label{lemma 5.2}
Let $\eta_0$ be distributed with $\nu_\gamma$. For any given $T>0$, as $N\rightarrow+\infty$,
\[
\left\{\frac{1}{\Lambda_{2, \alpha}(N)}M_{tN}^N:~0\leq t\leq T\right\}
\]
converges weakly, with respect to the Skorohod topology, to
\[
\left\{\sqrt{2f_1^{2, \alpha}(\mathbf{0})\beta(\gamma)}\mathcal{W}_t:~0\leq t\leq T\right\}.
\]
\end{lemma}
\proof
According to an argument similar to that given in the proof of Lemma \ref{lemma 4.2}, we only need to check an analogue of \eqref{equ 4.9}, i.e.,
\begin{equation}\label{equ 5.3}
\lim_{N\rightarrow+\infty}\frac{N}{\Lambda_{2, \alpha}^2(N)}\sum_{x\in \mathbb{Z}^1}\sum_{y\neq x}\|y-x\|_2^{-(d+\alpha)}\left(g_N(y)-g_N(x)\right)^2=2f_1^{2, \alpha}(\mathbf{0}).
\end{equation}
We only deal with the case of $\alpha=2$ since the case $\alpha>2$ follows from a similar argument. According to an argument similar to that given in the proof of Lemma \ref{lemma 4.2}, to check \eqref{equ 5.3} in the case of $\alpha=2$, we only need to show that
\begin{equation}\label{equ 5.4}
\lim_{N\rightarrow+\infty}\frac{1}{\log (\log N)}\int_0^{+\infty}e^{-s/N}p_s^{2, 2}(\mathbf{0}, \mathbf{0})ds=f_1^{2, 2}(\mathbf{0}).
\end{equation}
By Lemma \ref{lemma LCLT}, for any $\epsilon>0$, there exists $M_4=M_4(\epsilon)<+\infty$ such that
\[
(s\log s)p_s^{2, 2}(\mathbf{0}, \mathbf{0})\in \left((1-\epsilon)f_1^{2, 2}(\mathbf{0}), (1+\epsilon)f_1^{2, 2}(\mathbf{0})\right)
\]
when $s\geq M_4$. Then, according to the fact that
\[
\int_0^{M_4}e^{-s/N}ds+\int_{\epsilon N}^{+\infty}e^{-s/N}\frac{1}{s}ds=O(1)
\]
and
\[
\int_{M_4}^{\epsilon N}\frac{1}{s\log s}ds=\log (\log s)\Big|_{M_4}^{\epsilon N}=\log(\log N)(1+o(1)),
\]
Equation \eqref{equ 5.4} follows from an argument similar to that leading to \eqref{equ 4.14} and the proof is complete.
\qed

At last, we prove Theorem \ref{theorem 2.2 main theorem d=2} in the case $\alpha>2$.

\proof[Proof of Theorem \ref{theorem 2.2 main theorem d=2} in the case $\alpha\geq 2$]
Let $\Xi$ be the $(2, \alpha)$-version of that defined as in the proof of Theorem \ref{theorem 2.1 main theorem d=1} in the case $\alpha=1$, then we have an analogue of \eqref{equ 4.15}, i.e.,
\begin{equation}\label{equ 5.5}
\lim_{s\rightarrow+\infty}h^2_\alpha(s){\rm Var}_{\nu_\gamma}\left(\mathbb{E}_{\eta_0}\Xi(\eta_s)\right)=0.
\end{equation}
According to an argument similar to that given in the proof of Theorem \ref{theorem 2.1 main theorem d=1} in the case $\alpha=1$, Theorem \ref{theorem 2.2 main theorem d=2} in the case $\alpha\geq 2$ follows from Lemma \ref{lemma 5.2} and \eqref{equ 5.5}.
\qed

\quad

\textbf{Declaration of competing interest.}
The authors declare that they have no known competing financial interests or personal
relationships that could have appeared to influence the work reported in this paper.

\quad

\textbf{Data Availability.} Data sharing not applicable to this article as no datasets were generated or analysed during the current study.

\quad

\textbf{Acknowledgments.} The author is grateful to financial supports from the National Natural Science Foundation of China with grant number 12371142.

{}

\begin{thebibliography}{}

\bibitem{Bernardin2016}Bernardin, C., Gon\c{c}alves, P. and Sethuraman, S. (2016). Occupation times of long-range exclusion and connections to KPZ class exponents. \emph{Probability Theory and Related Fields} \textbf{166}, 365-428.

\bibitem{Dur2010}Durrett, R. (2010). Probability: Theory and Examples. 4th. Cambridge.

\bibitem{Ethier1986}Ethier, N. and Kurtz, T. (1986). \emph{Markov Processes: Characterization and Convergence.} John Wiley and Sons, Hoboken, NJ, USA.

\bibitem{Janvresse1999}Janvresse, E., Landim, C., Quastel, J. and Yau, H. (1999). Relaxation to equilibrium of conservative dynamics \uppercase\expandafter{\romannumeral1}: zero-range processes. \emph{The Annals of Probability} \textbf{27}, 325-360.

\bibitem{Kallenberg1997}Kallenberg, O. (1997). \emph{Foundations of Modern Probability.} Springer, New York.

\bibitem{Kipnis1987}Kipnis, C. (1987). Fluctuations des temps d'occupation d'un site dans l'exclusion simple sym\'{e}trique. \emph{Annales de l'IHP Probabilit\'{e}s et statistiques} \textbf{23}(1), 21-35.

\bibitem{kipnis+landim99} Kipnis, C. and Landim, C. (1999) \emph{Scaling limits of interacting particle systems.} Springer-Verlag, Berlin.

\bibitem{Kipnis1986}Kipnis, C. and Varadhan, S.R.S. (1986). Central limit theorem for additive functionals of reversible markov processes. \emph{Communications in Mathematical Physics} \textbf{104}, 1-19.

\bibitem{Komorowski2012}Komorowski, T., Landim, C. and Olla, S. (2012). \emph{Fluctuations in Markov processes: Time symmetry and Martingale approximation.} Springer, New York.

\bibitem{Lawler2010} Lawler, G., and Limic, V. (2010). \emph{Random Walk: a modern introduction.} Cambridge.

\bibitem{Quastel2002}Quastel, J.,  Jankowski, H. and Sheriff, J. (2002). Central limit theorem for zero-range processes. \emph{Methods and Applications of Analysis} \textbf{9}, 393-406.

\bibitem{Xue2025}Xue, X. (2025). Stationary fluctuations for occupation times of the long-range voter models on lattices. Arxiv: 2509.17518.

\bibitem{Zhao2025}Zhao, L. (2025). Stationary fluctuations for a multi-species zero range process with long jumps. \emph{Stochastic Processes and their Applications} \textbf{186}, 104645.


\end{thebibliography}
\end{document}